\documentclass[11pt]{article}
\usepackage[dvips]{graphics}
\usepackage{epsfig}%,rotating}
\usepackage{amssymb,amsmath}
\usepackage{latexsym}
\usepackage[usenames]{color}

\numberwithin{equation}{section}

\newcommand{\be}{\begin{equation}}
\newcommand{\ee}{\end{equation}}
\newcommand{\beaa}{\begin{eqnarray*}}
\newcommand{\eeaa}{\end{eqnarray*}}
\newcommand{\bea}{\begin{eqnarray}}
\newcommand{\eea}{\end{eqnarray}}
\newcommand{\lbl}{\label}

\newcommand{\bei}{\begin{itemize}}
\newcommand{\eei}{\end{itemize}}

\newcommand{\ml}{\mathcal}

\newcommand{\bd}{\bold}

\newtheorem{theorem}{ \noindent T{\footnotesize HEOREM}}
\newtheorem{prop}{ \noindent P{\footnotesize ROPOSITION}}
\newtheorem{lemma}{ \noindent L{\footnotesize EMMA}}[section]
\newtheorem{coro}{ \noindent C{\footnotesize OROLLARY}}
\begin{document}

%Coherence-Random-Matrix02092011
%\title{On the Coherence of High-Dimensional Random Matrices}
\title{Approximation of Rectangular Beta-Laguerre Ensembles and Large Deviations}
\author{Tiefeng Jiang$^{1}$ and Danning Li$^{2}$\\
University of Minnesota}

\date{}
\maketitle

\footnotetext[1]{ School of Statistics, University of Minnesota, 224 Church
Street S. E., MN55455, jiang040@umn.edu. \newline  \indent \ \
 The research of Tiefeng Jiang was
supported in part by NSF FRG Grant DMS-0449365 and NSF Grant DMS-1208982.  }
\footnotetext[2]{School of Statistics, University of Minnesota, 224 Church
Street S. E., MN55455 and Statistical Laboratory
 Centre for Mathematical Sciences,
Wilberforce Road, Cambridge
CB3 0WB,
England.

\ \
lixx0700@umn.edu.}

\begin{abstract}
\noindent Let $\lambda_1, \cdots, \lambda_n$ be random eigenvalues coming from the beta-Laguerre ensemble with parameter $p$, which is a generalization of the real, complex and quaternion  Wishart matrices of parameter $(n,p).$
%where $p$ is the dimension of the population distribution and $n$ is the sample size.
In the case that the sample size $n$ is much smaller than the dimension of the population distribution $p$, a common situation in modern data, we approximate the  beta-Laguerre ensemble by a beta-Hermite ensemble which is a generalization of the real, complex  and quaternion  Wigner matrices.  As corollaries, when $n$ is much smaller than $p,$ we show that the largest and smallest eigenvalues of the complex  Wishart matrix  are asymptotically independent; we obtain the limiting distribution of the condition numbers as a sum of two i.i.d. random variables with a Tracy-Widom distribution, which is much different from the exact square case that $n=p$ by Edelman (1988); we propose a test procedure for a spherical hypothesis test.  By the same approximation tool, we obtain the asymptotic distribution of the smallest eigenvalue of the beta-Laguerre ensemble.

In the second part of the paper, under the assumption that $n$ is much smaller than $p$ in a certain scale, we  prove the large deviation principles for three basic statistics: the largest eigenvalue, the smallest eigenvalue and the empirical distribution of $\lambda_1, \cdots, \lambda_n$, where the last large deviation is derived by using a non-standard method.
\end{abstract}

\bigskip



\newpage
\section{Introduction}
\lbl{intro}
\setcounter{equation}{0}
%\subsection{Background}

With the development of modern technology,  high-dimensional datasets appear very frequently in different scientific disciplines such as climate studies, financial data, information retrieval/search engines and functional data analysis.  The corresponding statistical problems have the feature that the dimension $p$ is possibly larger than the sample size $n.$ In particular, such feature  is very common in the data of gene expression. For example, in the data of ``1000 Genomes Project"  which is  by far the most detailed catalogue of human genetic variation, $n$ is usually at the level of  $10^3$ and $p$ is at the level of $10^7$ or $10^8$ (Durbin et al (2010)). In such cases the classical  statistical procedures based on fixed $p$ and large $n$ are no more applicable. The applications thus request new theories.  For recent progress in this area, see, for example, Candes and Tao (2005), Donoho et al (2006), Cai and Jiang (2011, 2012), and Vershynin (2012).

In this paper, we study the spectral properties of a Wishart matrix formed by a random sample of $p$-dimensional data with sample size $n$, where $p$ is larger than $n.$ Wishart matrices are very popular and useful objects in multivariate analysis, see, for example, the classical books by  Muirhead (1982) and Anderson (1984). It usually comes from the following formulation in statistics.
Let $\bd{y}_1, \cdots, \bd{y}_m$ be i.i.d. random variables with the $p$-dimensional multivariate normal distribution $N_p(\bd{\mu}, \bd{I}_p).$ Recall the sample covariance matrix
\begin{eqnarray}\lbl{coffee3}
 \bd{S}=\frac{1}{m}\sum_{i=1}^m(\bd{y}_i- \bar{\bd{y}})(\bd{y}_i- \bar{\bd{y}})^*\ \ \mbox{where}\ \  \bar{\bd{y}}=\frac{1}{m}\sum_{i=1}^m\bd{y}_i.
\end{eqnarray}
Then $m\bd{S}$ and $\bd{W}:=\bd{X}^*\bd{X}$ have the same distribution, where $\bd{X}=(x_{ij})_{n\times p}$, the random variables  $x_{ij}$'s are i.i.d. with distribution $N(0, 1)$ and $n=m-1.$ The matrix $\bd{W}$ is referred to as  the real Wishart matrix ($\beta=1$). If $x_{ij}$ are i.i.d. with the standard complex   or quaternion normal  distribution, then $\bd{W}$ is a complex Wishart matrix ($\beta=2$) or quaternion Wishart matrix ($\beta=4$).
%Let $\bd{X}=(x_{ij})$ be an $n\times p$ matrix with $n\geq p$, where $x_{ij}$'s are independent and identically distributed random variables with %the standard real normal ($\beta=1$), complex normal ($\beta=2$) or quaternion normal ($\beta=4$) distribution. Then $S_n=\bd{X}^*\bd{X}$ is a %$p\times p$ Wishart matrix.

Assume $p>n.$ Let $\lambda_1> \cdots> \lambda_n>0$ be the positive eigenvalues of  $\bd{W},$ which are the same as the $n$ eigenvalues of $\bd{X}\bd{X}^*.$ Write $\lambda=(\lambda_1, \cdots, \lambda_n).$ It is known that the density function of $\lambda$ is given by
\begin{eqnarray}\lbl{bWishart}
f_{n,\beta}(\lambda)=c_{n}^{\beta, p}\prod_{1\leq i<j\leq
n}|\lambda_i-\lambda_j|^{\beta}\cdot\prod_{i=1}^n\lambda_i^{\frac{\beta}{2}(p-n+1)-1}\cdot
e^{-\frac{1}{2}\sum_{i=1}^n\lambda_i}
\end{eqnarray}
for all $\lambda_1>0, \cdots, \lambda_n>0,$ where
\begin{eqnarray}\lbl{conWishart}
c_{n}^{\beta,
p}=2^{-\frac{\beta}{2} np}\prod_{j=1}^{n}\frac{\Gamma(1+\frac{\beta}{2})}
{\Gamma(1+\frac{\beta}{2}j)\Gamma(\frac{\beta}{2}(p-n+j))}.
\end{eqnarray}
See, for example,  James (1964) and Muirhead (1982)  for the cases $\beta=1$ and $2$,
and Macdonald (1995) and Edelman and Rao (2005) for $\beta=4$.
%, or (4.5) and (4.6) from \cite{Edelman05}.
The function $f_{n,\beta}(\lambda)$ in (\ref{bWishart}), being a probability density function for any $\beta>0,$  is called the $\beta$-Laguerre ensemble in literature. See, for example, Dumitriu (2003) and
Dumitriu and Edelman (2006).
% and Ram\'{i}rez,  Rider and  Vir\'{a}g (2011).

In this paper we will study the properties of $\lambda=(\lambda_1, \cdots, \lambda_n)$ for all $\beta>0.$ Precisely, there are two objectives. First, we show in Theorem \ref{main} that, when $p$ is much larger than $n$ in a certain scale, a ``normalized" $\beta$-Laguerre ensemble can be roughly thought as a $\beta$-Hermite ensemble with density function
\bea\lbl{betaHermite}
f_{\beta}(\lambda)=K_n^{\beta} \prod_{1\leq i <  j \leq n}|\lambda_i-\lambda_j|^{\beta}\cdot e^{-\frac{1}{2}\sum_{i=1}^n\lambda_i^2}
\eea
for all $\lambda=(\lambda_1, \cdots, \lambda_n)\in \mathbb{R}^n,$ where
\bea\lbl{Hermitecons}
K_n^{\beta}=(2\pi)^{-n/2}\prod_{j=1}^n\frac{\Gamma(1+\frac{\beta}{2})}{\Gamma(1+\frac{\beta}{2}j)}.
\eea
The eigenvalues $\lambda_1,\cdots, \lambda_n$ of
the Gaussian orthogonal ensemble (GOE), Gaussian unitary ensemble (GUE) and Gaussian symplectic ensemble (GSE) have the joint density function $f_{\beta}(\lambda)$ as in (\ref{betaHermite})  with $\beta=1, 2$ and $4,$ respectively. See, e.g., chapter 17 from Mehta (1991) for more details.

By using Theorem \ref{main} mentioned above (\ref{betaHermite})
 and some known results on $\beta$-Hermite ensembles, under the assumption that $p$ is much larger  than $n$ in a certain scale, we obtain the following new results:

(i) The largest eigenvalue $\lambda_{max}$ and smallest eigenvalue $\lambda_{min}$ in the $\beta$-Laguerre ensemble as $\beta=2$ are asymptotically  independent (Proposition \ref{independence}).

(ii) The condition number $\kappa_n$ (see the definition in (\ref{Changchun})) of an $n\times p$ matrix $(x_{ij})$ ($x_{ij}$'s are i.i.d. centered complex Gaussian random variables), when suitably normalized, converges weakly to $U+V$ where $U$ and $V$ are independent random variables with a common Tracy-Widom law  (Corollary \ref{prasing}). This is much different from the exact square case that $n=p$ studied by Edelman (1988):  $\kappa_n/n$ converges weakly to a distribution with density function $h(x)=8x^{-3}e^{-4/x^2}$ for $x>0.$ Based on this result, a spherical test in statistics is proposed below Corollary \ref{prasing}.

(iii) A linear transform of the smallest eigenvalue  of the $\beta$-Laguerre ensemble converges to the $\beta$-Tracy-Widom law for any $\beta>0$ (Proposition \ref{singing}). The counterpart for the largest eigenvalues was studied by Ram\'{i}rez et al (2011).

It is worthwhile to state that the condition number $\kappa_n$ mentioned in (ii) is an important quantity in the field of numerical analysis dated back to Von Neumann and Goldstine (1963).

%According to Bai and Yin (1988), the empirical distribution of the eigenvalues of $\sqrt{\frac{n}{p}}(\frac{S_n}{n}-I_n)$ converges weakly to the %semicircle law with pdf
%\beaa
%f(x)=\frac{1}{2\pi}\sqrt{4-x^2},\ \, |x|\leq 2.
%\eeaa

In the second part of this paper we study the large deviations for the eigenvalues of the $\beta$-Laguerre ensembles when $p$ is much larger than $n.$ The large deviation for the eigenvalues of random matrices is one of active research areas in random matrix theory. See, for example, a survey paper by Guionnet (2004) and some chapters from  Hiai and Petz (2006) and Anderson et al (2009). In particular, Ben Arous and Guionnet (1997) investigate the large deviation for Wigner's semi-circle law; Ben Arous,  Dembo and Guionnet (2001) and Anderson et al (2009) study the largest eigenvalues of Wigner and Wishart matrices;   A corollary from Hiai  and Petz (1998) says that a normalized empirical distribution $\mu_n$ of the positive eigenvalues of real Wishart matrix $\bd{X}_{n\times p}\bd{X}^*_{n\times p}$ follows the large deviation principle (LDP) such that
\bea
& & \limsup_{n\to\infty}\frac{1}{p^2}\log P(\mu_n \in F) \leq -\inf_{\nu\in F}I(\nu) \ \ \mbox{and} \lbl{upperb}\\
& & \liminf_{n\to\infty}\frac{1}{p^2}\log P(\mu_n \in G) \geq -\inf_{\nu\in G}I(\nu) \lbl{lowerb}
\eea
for every closed set $F$ and open set $G$ under the topology of weak convergence of probability measures on $\mathbb{R},$ where $\mu_n:=\frac{1}{n}\sum_{i=1}^n\delta_{\lambda_i/p}$, the eigenvalues $\lambda_1, \cdots, \lambda_n$ have joint density $f_{n,1}(\lambda)$ as in (\ref{bWishart}) with $n/p\to \gamma\in (0,1],$ and
\beaa
I(\nu)=-\frac{\gamma^2}{2}\iint\log |x-y|\, d\nu(x)d\nu(y) + \frac{\gamma}{2}\int\left(x-(1-\gamma)\log x\right)d\nu(x) + \mbox{Const}
\eeaa
which takes the minimum value zero at  the Marchenko-Pastur law  with  density function
\bea\lbl{robe}
h(x)=\frac{1}{2\pi \gamma x}\sqrt{(x-\gamma_1)(\gamma_{2}-x)}
\eea
for $x\in [\gamma_{1}, \gamma_{2}]$ and $\gamma_1=(\sqrt{\gamma}-1)^2$ and $\gamma_2=(\sqrt{\gamma}+1)^2.$  For the general framework of the large deviation principle, its connection to the subjects of mathematics, physics, statistics and engineering,  see, for example, Shwartz and Weiss (1995),  Dembo and Zeitouni (2009) and Ellis (2011).

When $p/n\to \infty$, the LDP problem for $\mu_n$ in (\ref{upperb}) and (\ref{lowerb}) has been open until now. In fact, we resolve the problem in Theorem \ref{main1}  under the assumption that both $p$ and $n$ are large with $p/n^2\to \infty.$ Contrary to the Marchenko-Pastur law stated in (\ref{robe}), we show that the rate function in Theorem \ref{main1} takes the minimum value at the semi-circle law.

The large deviation principles for the largest eigenvalue $\lambda_{max}$ and the smallest eigenvalue $\lambda_{min}$ of the $\beta$-Laguerre ensemble are also studied in Theorems \ref{pack} and \ref{done1} as $p/n\to \infty.$ Their rate functions are explicit.

The rest of the paper is organized as follows. In Section \ref{convergence} we give a theorem that the $\beta$-Laguerre ensemble converges to the $\beta$-Hermite ensemble as $p$ is much larger than $n$ and present some implications; In section \ref{LDPsection} we give three  theorems about the large deviations for the largest eigenvalues, the smallest eigenvalues and  the empirical distributions of the eigenvalues of the $\beta$-Laguerre ensemble as $p$ is much larger than $n$ in a certain scale. The proofs of the results stated in Sections \ref{convergence} and \ref{LDPsection} are given in Sections \ref{Proofconvergence} and \ref{Proofrainbow}, respectively.

The reader is warned that the notation $\mu$ or $\mu_n$ throughout the paper sometimes represents a probability measure, a mean value or an eigenvalue in different occasions, but this will not cause confusions from the context.

\section{Convergence of Laguerre Ensembles to Hermite Ensembles and Its Applications}\lbl{convergence}

Let $\mu$ and $\nu$ be probability measures on $(\mathbb{R}^k, \mathcal{B}),$ where $k\geq 1$ and $\mathcal{B}$ is the Borel $\sigma$-algebra on $\mathbb{R}^k$. The variation distance $\|\mu - \nu\|$ is defined by
\begin{eqnarray}\lbl{variation}
\|\mu- \nu\|= 2\cdot \sup_{A \in \mathcal{B}}|\mu(A) - \nu (A)|= \int_{\mathbb{R}^k}|f(x)-g(x)|\,dx_1\cdots dx_k
\end{eqnarray}
if $\mu$ and $\nu$ have density functions $f(x)$ and $g(x)$ with respect to the Lebesgue measure. For a random vector $Z$, we use $\mathcal{L}(Z)$ to denote its probability distribution. The notation $a_n \gg b_n$ means $\lim_{n\to\infty}a_n/b_n=+\infty.$

\begin{theorem}\lbl{main} Let $\lambda=(\lambda_1, \cdots, \lambda_{n})$ be random variables with density function $f_{n,\beta}(\lambda)$ as in (\ref{bWishart})  and $\mu=(\mu_1, \cdots,
\mu_{n})$ be random variables with density $f_{\beta}(\mu)$ as in (\ref{betaHermite}). Set $x_i=\sqrt{\frac{p}{2\beta }}\big(\frac{\lambda_i}{p}-\beta\big)$ for $1\leq i \leq n.$ Then $\|\ml{L}\big((x_1, \cdots, x_n)\big)-\ml{L}\big((\mu_1, \cdots, \mu_n)\big)\| \to 0$ if (i) $n\to\infty$ and $p=p_n \gg n^3$ or (ii) $n$ is fixed and $p\to\infty.$
\end{theorem}

\begin{figure}[t]
\centerline{\includegraphics[scale=0.77]{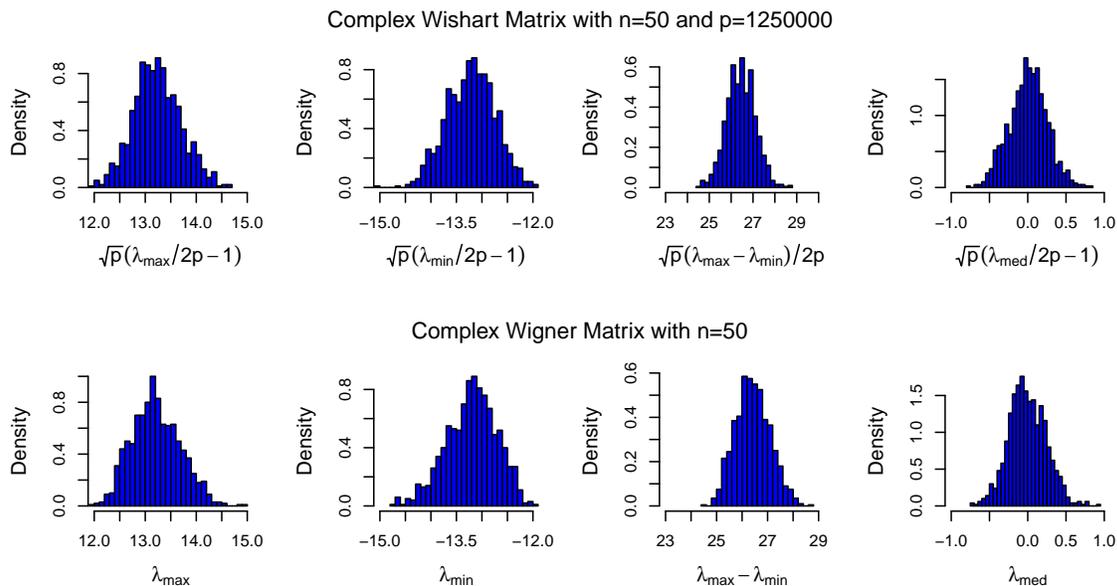}}
\caption{\sl Histogram  of Theorem \ref{main}. Simulation is run to compare distribution of eigenvalue of Laguerre ensemble with $n=50$, $p=1250000$ and  $\beta=2$ and those of Hermite ensemble with $n=50$. Picture on top row show empirical distribution of largest, smallest, range, and median of eigenvalue of Laguerre ensemble with suitable normalization. Picture on bottom row show corresponding distribution of Hermite ensemble.}
%Picture on top row show empirical distribution of largest, smallest, range and  median of eigenvalue  with  suitable normalization. Pictures on bottom row show %limit distribution of top row. It is seen, when $p$ is much larger than $n$, they  behave similarly.}
\label{fig1}
\end{figure}

Roughly speaking, Theorem \ref{main} says that a ``very rectangular-shaped" $\beta$-Laguerre ensemble is essentially a $\beta$-Hermite ensemble. This can be seen from a simulation study as shown in Figure \ref{fig1}. It indicates that the maximum, the minimum, the median and the range of the former one are close to those of the latter one. Theoretically, comparing the density functions of both ensembles, with the transform given in Theorem \ref{main}, the essential understanding is that $(1+\frac{\lambda}{n})^ne^{-\lambda} \sim e^{-\lambda^2/(2n)}$ as $n\to\infty.$ The left hand side is from the pdf of the Laguerre ensemble, and the right is from the Hermite ensemble. The term $J:=\prod_{1\leq i<j\leq n}|\lambda_i-\lambda_j|^{\beta}$ under the linear transform $a\lambda+b$ is identical to $a^{\beta n(n-1)/2}J.$ Multiplying both sides of ``$\sim$" by $J$, we then see that the two density functions are asymptotically identically as $p$ is much larger than $n$. The literal argument is given in (\ref{shaq}).

 Theorem \ref{main} can also be understood through random matrix models. For example, take the Wishart matrix $\bd{W}=\bd{X}\bd{X}^*$ where $\bd{X}=(x_{ij})_{n\times p}$ and $x_{ij}$'s are i.i.d. random variables with the distribution of $N(0, 1)$ and $p>n$. For simplicity, let $n$ be fixed, then the entries of $(\bd{W}-p\bd{I})/\sqrt{p}$ are asymptotically independent normals. In other words, the limit is an Hermite ensemble. Finally, according to Dyson's "three fold way" of classical random matrices (Dyson, 1962), three types of random matrices are of great interest: Hermite, Laguerre and Jacobi matries. This paper together with the work by Jiang (2009, 2013) says that, as the parameters of the last two matrices are in extreme relationships, we have
\beaa
\mbox{Jacobi}\ \to\ \mbox{Laguerre}\ \to \ \mbox{Hermite}.
\eeaa
See Theorem \ref{main} here and Theorem 6 from Jiang (2013). The notation ``$\to$" is interpreted by ``is reduced to" in words or the variation approximation as in Theorem \ref{main} literally.

Now, by combining Theorem \ref{main} with some known results on the $\beta$-Hermite ensembles, we obtain several new results. To state them, let us first review the Tracy-Widom distributions. Set

\bea\lbl{F2}
F_2(x)=\exp\Big(-\int_x^{\infty}(y-x)q^2(y)\,dy\Big),\  x\in \mathbb{R},
\eea
where $q$ is the unique solution to the Painlev\'{e} II differential equation
\bea\lbl{Pinaleve}
q''(x)=xq(x) + 2q^3(x)
\eea
satisfying the boundary condition $q(x)\sim Ai(x)$ as $x\to\infty,$ where $Ai(x)$ is the Airy function.  It is known from Hastings and McLeod (1980) that
\beaa
q(x)=\sqrt{-\frac{x}{2}}\Big(1+\frac{1}{8x^3} + O\big(\frac{1}{x^6}\big)\Big)
\eeaa
as $x\to -\infty.$ The distributions for the orthogonal and symplectic cases (Tracy and Widom (1996)) are
\bea
& & F_1(x)=\exp\Big(-\frac{1}{2}\int_x^{\infty}q(y)\,dy\Big)(F_2(x))^{1/2}\ \ \mbox{and}\  \lbl{F1} \\
& & F_4(x/\sqrt{2})=\cosh\Big(\frac{1}{2}\int_x^{\infty}q(y)\,dy\Big)(F_2(x))^{1/2}\lbl{F4}
\eea
for all $x\in \mathbb{R},$ where $\cosh t=(e^{t}+ e^{-t})/2$ for $t\in \mathbb{R}.$

Our first result following from Theorem \ref{main} is on complex Wishart matrices ($\beta=2$). In fact, data matrices with complex-valued
entries arise frequently, for example, in signal processing applications (e.g.,
 Couillet and Debbah, 2011) and statistics (e.g., James (1964) and Picinbono (1996)). Given $\lambda_1, \cdots, \lambda_n, $ set $\lambda_{min}=\min\{\lambda_1, \cdots, \lambda_n\}$ and   $\lambda_{max}=\max\{\lambda_1, \cdots, \lambda_n\}.$

\begin{prop}\lbl{independence} Let $\lambda=(\lambda_1, \cdots, \lambda_{n})$ be random variables with density  $f_{n,\beta}(\lambda)$ as in (\ref{bWishart}) with $\beta=2$. Set $\mu_{n,1}=2p-4\sqrt{np}$, $\mu_{n,2}=2p+4\sqrt{np}$ and $\sigma_n=2\sqrt{p}n^{-1/6}.$ If $p=p_n\to\infty$ and $p \gg n^3,$ then $\big((\lambda_{min}-\mu_{n,1})/\sigma_n, (\lambda_{max}-\mu_{n,2})/\sigma_n\big) \in \mathbb{R}^2$  converge weakly to $(-U, V)$,  where  $U$ and $V$ are i.i.d. with distribution function $F_2(x)$ as in (\ref{F2}).
\end{prop}

Basor et al (2012) heuristically show that Proposition \ref{independence} is true as $n/p \to c$. Our result above is rigorous. It remains open  at this moment if Proposition \ref{independence} is still true for $\beta \ne 2$ under the assumption  $p/n\to \gamma \in [1, \infty].$

Let $\bd{X}=(x_{ij})_{n\times p}$ and $x_{ij}$'s be i.i.d. complex random variables with the distribution of $(\xi+i\eta)/\sqrt{2}$ where $\xi$ and $\eta$ are i.i.d. with $N(0,1)$-distribution. Suppose $\lambda_{max}$ and $\lambda_{min}$ are the largest and smallest eigenvalues of the matrix $\bd{X}\bd{X}^*$ (the density of the eigenvalues of this matrix corresponds to $\beta=2$ in (\ref{bWishart})).
%An immediate consequence of Proposition \ref{independence} is the following result about the condition number %$\kappa_n:=(\lambda_{max}/\lambda_{min})^{1/2}.$
Recall the condition number defined by
\bea\lbl{Changchun}
\kappa_n:=\Big(\frac{\lambda_{max}}{\lambda_{min}}\Big)^{1/2}.
\eea
An immediate consequence of Proposition \ref{independence} is the following result about $\kappa_n$.
\begin{coro}\lbl{prasing} Let $\lambda=(\lambda_1, \cdots, \lambda_{n})$ be the eigenvalues of $\bd{X}\bd{X}^*$ stated above. If $p\gg n^3,$
%random variables with density  $f_{n,\beta}(\lambda)$ as in (\ref{bWishart})  with $\beta=2$ .
then $\alpha_n(\kappa_n -\beta_n)$ converges weakly to $U+V,$ where $\alpha_n=2\sqrt{p}n^{1/6}$, $\beta_n=1+2\sqrt{\frac{n}{p}}$, and $U$ and $V$ are i.i.d. with distribution function $F_2(x)$ as in (\ref{F2}).
\end{coro}

In the exact square case that $p=n$, Edelman (1988) proves that $\kappa_n/n$ converges weakly to a distribution with density function $h(x)=8x^{-3}e^{-4/x^2}$ for $x>0.$ In the rectangular case such that $p\gg n^3,$ Corollary \ref{prasing}  shows a very different behavior of $\kappa_n$.

\begin{figure}
%\centerline{
\includegraphics[scale=0.77]{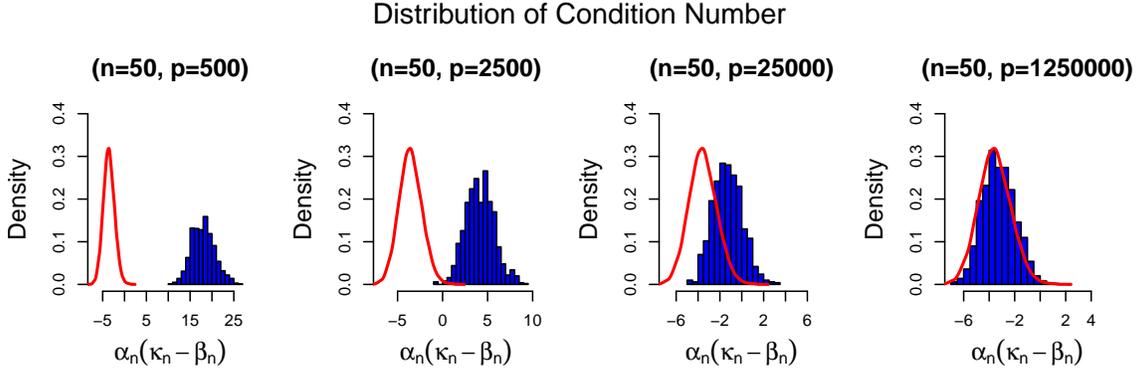}%}
\caption{\sl Histograms of Corollary \ref{prasing}. Simulation is run with $n=50$ and $p=500, 2500, 25000, 125000.$ Solid curve, as density  of  $U+V$,  fit well to histogram of $ \alpha_n(\kappa_n-\beta_n)$  when $p$ become large. This is consistent with Corollary \ref{prasing}. The solid curve is simulated according to the R-package by Johnstone et al, see the link in the reference.}
\label{fig2}
\end{figure}

Let $\bd{Y}=\bd{\xi}+i\bd{\eta}$ be a multivariate complex normal distribution where $\bd{\xi}\sim N_p(\bd{0}, \bd{\Sigma}_1)$ and $\bd{\eta}\sim N_p(\bd{0}, \bd{\Sigma}_2)$ are independent. Consider the spherical test $H_0: \bd{\Sigma}_1=\bd{\Sigma}_2=\rho \bd{I}_p$ vs $H_a: H_0$ is not true, where $\rho>0$ is not specified. Let $\bd{Y}_1, \cdots, \bd{Y}_n$ be a random sample from the population distribution of $\bd{Y}$ with $p>n.$ The classical likelihood ratio test (see, e.g., Muirhead (1982)) does not work here simply because  $p>n$. Ledoit and Wolf (2002) suggest to use certain alternative test statistics when $p/n\to y\in (0, \infty)$. Jiang and Yang (2013) and Jiang and Qi (2013) investigate the test through obtaining the central limit theorem of its likelihood ratio test statistic under $p/n\to y\in [0, 1]$. Dette and Holger  (2005) further extend Ledoit and Wolf's result to cover the case of $y =0$ or $\infty$.  Corollary \ref{prasing} in our paper supplements these results by providing another way of running the sphericity test when $p$ is much larger than $n$ and the data are complex-valued. In fact, we can carry the test in the following way. Set $\bd{X}=\bd{X}_{n\times p}=(\bd{Y}_1, \cdots, \bd{Y}_n)'.$ Then the $n$ positive eigenvalues of $\bd{X}^*\bd{X}/\rho^2$ have the joint density function $f_{n,\beta}(\lambda)$  as in (\ref{bWishart})  with $\beta=2$. Recall that $\lambda_{max}$ and $\lambda_{min}$ are the largest and smallest positive eigenvalues of $\bd{X}^*\bd{X},$ respectively. Then, $\kappa_n=(\lambda_{max}/\lambda_{min})^{1/2}$ does not depend on the unknown parameter $\rho.$  By Corollary \ref{prasing}, the region to reject $H_0$ with an asymptotic $1-\alpha$ confidence level is $\{\alpha_n|\kappa_n -\beta_n|>s\}$, where $s>0$ satisfies $P(|U+V|>s)=\alpha.$ The value of $s$ can be calculated through a numerical method by using (\ref{F2}), (\ref{Pinaleve}) and the independence between $U$ and $V$.

%\begin{prop}\lbl{dissertation} Let $\lambda=(\lambda_1, \cdots, \lambda_{p})$ be random variables with density  $f_{n,\beta}(\lambda)$ as in %(\ref{bWishart}). Set $\mu_n=\beta (n-2 \sqrt{np})$ and $\sigma_n=\beta \sqrt{n}p^{-1/6}.$ If $p=p_n\to\infty$ and $p =o(n^{1/3}),$ then, as %$n\to\infty$, $(\lambda_{min}-\mu_n)/\sigma_n$
%%\beaa
%%\frac{p^{1/6}}{\beta \sqrt{n}}\big(\beta n-2\beta p\sqrt{n}-\lambda_{min}\,\big)
%%\eeaa
%converges weakly to the distribution function $1-F_{\beta}(-x),\, x\in \mathbb{R},$ where $F_{\beta}(x)$ is as in (\ref{F1}), (\ref{F2}) and %(\ref{F4})  for $\beta=1,2$ and $4.$
%\end{prop}
%Assuming  $\beta=1$ and $n/p\to \gamma \in (1, \infty),$ Ma (2010) shows that  $(\log \lambda_{min}-\nu_n)/\tau_n$ converges weakly to the %distribution function $1-F_1(-x)$ as $n\to\infty,$ where $\nu_n$ and $\tau_n$ are normalizing constants.  Proposition \ref{dissertation} gives %the limiting distribution of $\lambda_{min}$ when $\gamma=\infty$ and $\beta=1, 2, 4.$ For the largest eigenvalue $\lambda_{max},$ Johansson %(2000), Johnstone (2001) and Karoui (2003) obtain its limiting distribution as $\beta=1,2, 4$ and $\gamma\in [0, \infty].$

\begin{figure}[t]
\centerline{\includegraphics[scale=0.77]{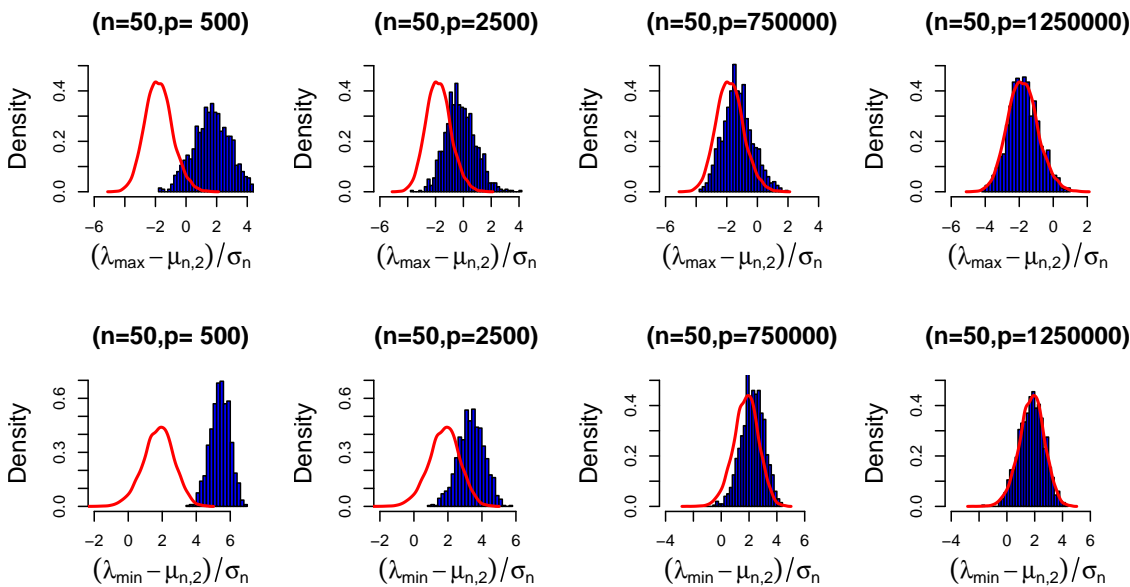}}
\caption{\sl Histogram of Proposition \ref{singing} and discussion below. Simulation is run to compare Tracy-Widom law (solid curve) and histogram of largest and smallest eigenvalue with $\beta=2$. It is seen that Tracy-Widom law fit histogram well with increase of $p$.  }
\label{fig3}
\end{figure}
Now we study the limiting distribution of $\lambda_{min}$ in the $\beta$-Laguerre ensemble for all $\beta >0.$ To do so, consider the random operator
\begin{eqnarray}\lbl{operator}
\mathcal{H}_{\beta}=-\frac{d^2}{dx^2} + x +\frac{2}{\sqrt{\beta}}b_x'
\end{eqnarray}
where $b_x$ is a standard Brownian motion on $[0, +\infty)$ ($b_x'$ is not the derivative of $b_x$ since it is not differentiable almost everywhere). We use equation (\ref{operator}) in the following sense. For $\lambda\in \mathbb{R}$ and  function $\psi(x)$ defined on $[0, +\infty)$ with $\psi(0)=0$ and $\int_0^{\infty}\left((\psi')^2 + (1+x)\psi^2\right)\,dx<\infty,$ we say $(\psi, \lambda)$ is an eigenfunction/eigenvalue pair for $\mathcal{H}_{\beta}$ if $\int_0^{\infty}\psi^2(x)\,dx=1$ and
\begin{eqnarray*}
\psi''(x)=\frac{2}{\sqrt{\beta}}\psi(x)b_x'+(x- \lambda)\psi(x)
\end{eqnarray*}
holds in the sense of integration-by-parts, that is,
\begin{eqnarray}\lbl{radish}
\psi'(x)-\psi'(0)=\frac{2}{\sqrt{\beta}}\psi(x)b_x+\int_0^x-\frac{2}{\sqrt{\beta}}b_y\psi'(y)\,dy + \int_0^x(y -\lambda)\psi(y)\,dy.
\end{eqnarray}
Theorem 1.1 from Ram\'{i}rez et al (2011) shows that, with probability one, for each $k\geq 1,$ the set of eigenvalues of $\mathcal{H}_{\beta}$ has well-defined $k$-lowest eigenvalues $(\Lambda_0, \cdots, \Lambda_{k-1}).$ Our result on $\lambda_{min}$ of a $\beta$-Laguerre ensemble is given next.

\begin{prop}\lbl{singing} Let $\beta>0$ and $\lambda=(\lambda_1, \cdots, \lambda_{n})$ be random variables with density  $f_{n,\beta}(\lambda)$ as in (\ref{bWishart}). Set $\mu_n=\beta (p-2 \sqrt{np})$ and $\sigma_n=\beta \sqrt{p}n^{-1/6}.$ If $p=p_n\to\infty$ and $p \gg n^3,$ then, as $n\to\infty$, $(\lambda_{min}-\mu_n)/\sigma_n$
%\beaa
%\frac{p^{1/6}}{\beta \sqrt{n}}\big(\beta n-2\beta p\sqrt{n}-\lambda_{min}\,\big)
%\eeaa
converges weakly to the distribution of $\Lambda_0.$
\end{prop}

%Indeed Proposition \ref{dissertation} is a special case of Proposition \ref{singing}. We differentiate the two cases because the former has %applications in Statistics for real, complex and quaternionic  Wishart matrices.
It is known from Ram\'{i}rez et al (2011) that $-\Lambda_0$ has the  distribution  $F_{\beta}(x)$  as in (\ref{F1}), (\ref{F2}) and (\ref{F4})  for $\beta=1,2$ and $4.$

Under the less restrictive condition that $p \gg n,$ Paul (2011) obtains Proposition  \ref{singing} for $\beta=1$ and $2.$ Further, assuming  $p/n\to \gamma \in (1, \infty),$ Baker et al (1998) show  that, if $\beta=2$, then $(\lambda_{min}-\nu_n)/\tau_n$ converges weakly to the distribution function $1-F_2(-x)$ (the distribution of $-\Lambda_0$ for $\beta=2$), where $\nu_n$ and $\tau_n$ are normalizing constants. Ma (2010) obtains a similar result for $\beta=1.$ Here, Proposition \ref{singing} holds for any $\beta>0.$
%Proposition \ref{dissertation} gives the limiting distribution of $\lambda_{min}$ when $\gamma=\infty$ and $\beta=1, 2, 4.$

For largest eigenvalue $\lambda_{max},$ Johansson (2000), Johnstone (2001) and Karoui (2003) obtain its limiting distribution as $\beta=1,2, 4$ and $\gamma\in [0, \infty].$ For general $\beta>0,$ the limiting distribution of $\lambda_{max}$ is obtained by
Ram\'{i}rez et al (2011)  for the $\beta$-Laguerre ensembles (that is, $\lambda=(\lambda_1, \cdots, \lambda_{n})$ has density  $f_{n,\beta}(\lambda)$ as in (\ref{bWishart})) as $p/n\to \gamma \in [1, \infty)$. We derive the asymptotic distribution of $\lambda_{min}$ for the same $\beta$-Laguerre ensemble when $p \gg n^3$ in Proposition \ref{singing}.  At this point it is not known if a result similar to Proposition \ref{singing} still holds as $p/n\to \gamma\in (0, \infty).$

Although the proof of Theorem \ref{main} suggests that the order of $p \gg n^3$ in Theorem \ref{main} is the best one to make the approximation hold, the orders appearing in Propositions \ref{independence} and \ref{singing} and  Corollary \ref{prasing} could be relaxed. This is because Theorem \ref{main} is a uniform approximation, and the three results are specific cases. One can see improvements in a different but similar situation in Dong et al (2012).

\section{Large Deviations for Eigenvalues}\lbl{LDPsection}
In this section we study the large deviations for three basic statistics as $p \gg n:$  the largest eigenvalue $\lambda_{max}$, the smallest eigenvalue $\lambda_{min}$ and the empirical distribution of $\lambda_1, \cdots, \lambda_n$ which come from a $\beta$-Laguerre ensemble.  One can check, for example, Dembo and Zeitouni (2009) for the definition of the large deviation principle (LDP). The first one is about the largest eigenvalue.

\begin{theorem}\lbl{pack} Suppose $\lambda_1, \cdots, \lambda_n$ have the density $f_{n,\beta}(\lambda)$ as in (\ref{bWishart}). Assume $p=p_n \gg n$ as $n\to\infty.$  Then, $\{\frac{\lambda_{max}}{p};\, n\geq 2\}$ satisfies the LDP with speed $\{p_n; \, n\geq 2\}$ and rate function $I(x)$ where
\beaa
I(x)=
\begin{cases}
\frac{x-\beta}{2} -\frac{\beta}{2}\log \frac{x}{\beta}, & \text{if $x\geq \beta$;}\\
+\infty, & \text{if $x<\beta$}.
\end{cases}
\eeaa
\end{theorem}
For the smallest eigenvalue, we have the following.
\begin{theorem}\lbl{done1} Suppose $\lambda_1, \cdots, \lambda_n$ have the density $f_{n,\beta}(\lambda)$ as in (\ref{bWishart}). Assume  $p=p_n \gg n$  as $n\to\infty.$ Then, $\{\frac{\lambda_{min}}{p};\, n\geq 2\}$ satisfies the LDP with speed $\{p_n; \, n\geq 2\}$ and rate function $I(x)$ where
\beaa
I(x)=
\begin{cases}
\frac{x-\beta}{2} -\frac{\beta}{2}\log \frac{x}{\beta}, & \text{if $0<x\leq \beta$;}\\
+\infty, & \text{otherwise}.
\end{cases}
\eeaa
\end{theorem}
\begin{figure}[t]
\centerline{\includegraphics[scale=0.77]{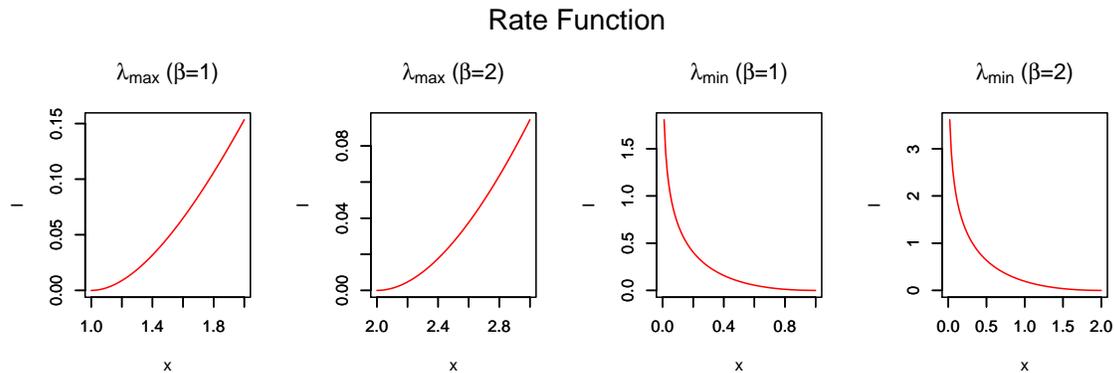}}
\caption{\sl  Rate function in Theorem \ref{pack} and \ref{done1} for largest and smallest eigenvalue with $\beta=1, 2$, respectively. }
\label{fig4}
\end{figure}
Since  $\lambda_{max}/p$ and $\lambda_{min}/p$ are all positive, it can be seen from Theorems \ref{pack} and \ref{done1} that the two rate functions on $(0, \infty)$ look ``symmetric" with respect to the line $x=\beta.$ When $p/n\to (0, \infty),$ the large deviations for $\lambda_{max}/p$ is known, see, for example, Anderson et al (2009). The LDP for $\lambda_{max}$   has a potential use in statistical tests. For example, Maida (2007)  investigates  the LDP for $\lambda_{max}$ of rank one deformations of Gaussian ensembles when $p/n$ converges to a constant (corresponding to the Laguerre ensemble with $\beta=2$). Bianchi et al (2011) uses Maida's work  to develop a statistical test for detecting single-source. One may use Theorem {\ref{pack}} and {\ref{done1}} in this paper to perform some statistical tests for the case $p\gg n$.

\noindent Now we consider the large deviation for the empirical distribution of the eigenvalues.

\begin{theorem}\lbl{main1} Given $\beta>0,$ let $\lambda=(\lambda_1, \cdots, \lambda_n)$ have the joint density function $f_{n,\beta}(\lambda)$ as in (\ref{bWishart}). Set $x_i=\sqrt{\frac{p}{2\beta }}\big(\frac{\lambda_i}{p}-\beta\big)$ for $1\leq i \leq n$  and
\bea\lbl{spring}
\mu_n=\frac{1}{n}\sum_{i=1}^n\delta_{\frac{x_i}{\sqrt{n}}}.
%\mu_n=\frac{1}{p}\sum_{i=1}^p\delta_{\sqrt{\frac{n}{p}}(\frac{\lambda_i}{n}-\beta)}.
\eea
If $p=p_n\to \infty$ and $p \gg n^2,$ then $\{\mu_n;\, n\geq 2\}$ satisfies the LDP with speed $\{n^2\}$ and rate function $I_{\beta}(\nu),$ where
\bea\lbl{LDPrate1}
I_{\beta}(\nu)=\frac{1}{2}\int_{\mathbb{R}^2}g(x,y)\,\nu(dx)\,\nu(dy) +\frac{\beta}{4}\log \frac{\beta}{2}-\frac{3}{8}\beta
\eea
and
\bea\lbl{elm1}
g(x,y)=
\begin{cases}
\frac{1}{2}(x^2+y^2) - \beta \log |x-y|, & \text{if $x\ne y$;}\\
+\infty, & \text{otherwise.}
\end{cases}
\eea
\end{theorem}

Theorem \ref{main1} is consistent with Theorem \ref{main} which says that the $\beta$-Laguerre ensemble is ``essentially" a $\beta$-Hermite ensemble when $p$ is large  enough relative to $n.$

The proof of Theorem \ref{main1} is different from the standard method of proving large deviations for the empirical distributions of eigenvalues (see, e.g., Ben Arous and Guionnet (1997), Hiai  and Petz (1998) and Guionnet (2004)). In fact, reviewing (\ref{upperb}) and (\ref{lowerb}), we estimate $P(\mu_n\in A)$ for a set $A$ by making a measure transformation such that the underlying $\beta$-Laguerre distribution is changed to a $\beta$-Hermite distribution. After an approximation step similar to that in Theorem \ref{main}, we use the known result on LDP for $\beta$-Hermite ensembles to complete the proof.

As stated in Theorem 1.3 from Ben Arous and Guionnet (1997), $I_{\beta}(\nu)$ is the rate function of the large deviation for $\mu_n$ in (\ref{spring}) with $p=n$ when  $x_1, \cdots, x_n$ come from a $\beta$-Hermite ensemble with density $f_{\beta}(x)$ as in (\ref{betaHermite}), and the rate function $I_{\beta}(\nu)$ takes the unique minimum value $0$ at the semi-circle law with density function $g_{\beta}(x)=(\beta \pi)^{-1}\sqrt{2\beta -x^2}$ for any $|x|\leq \sqrt{2\beta}$ and $\beta>0.$  This fact implies a weak convergence of the spectral distribution: under the setting of Theorem \ref{main1}, with probability one, $\mu_n$ converges weakly to a probability  distribution with density $g_{\beta}(x).$ When $\beta=1,$   the underlying matrix is the real Wishart matrix, Bai and Yin (1988) show the weak convergence  with the relaxed condition $p \gg n.$ Our next result says that for all $\beta >0$, if $p \gg n$, the limiting empirical distribution of the eigenvalues is still a semi-circle law with a different radius.

\begin{prop}\lbl{LSD}Given $\beta>0,$ suppose $\lambda=(\lambda_1, \cdots, \lambda_n)$ has the density  $f_{n,\beta}(\lambda)$ as in (\ref{bWishart}). Let $x_i $    and $\mu_n$ be as in (\ref{spring}).  If $p \gg n\to \infty$, then, with probability 1,  $\mu_n$ converges to the semi-circle law weakly with density  $g_{\beta}(x)=(\beta \pi)^{-1}\sqrt{2\beta -x^2}$ for $|x|\leq \sqrt{2\beta}$.
 \end{prop}

\begin{figure}[t]
\centerline{\includegraphics[scale=0.76]{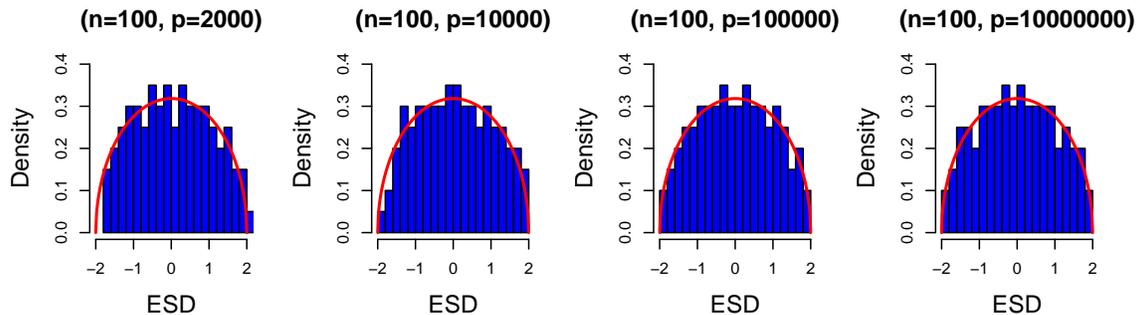}}
\caption{\sl Histogram of  Corollary \ref{LSD}. Picture  show that density of semicircle law (solid curve)  fit histogram of eigenvalue well even $p$ is not very  large, }
\label{fig5}
\end{figure}

Finally, the order $p \gg n$ in Theorems \ref{pack} and \ref{done1} is the best order to make the theorems hold when one considers the case $p$ being much larger  than $n.$ From the proof of Theorem \ref{main1}, we see the order $p \gg n^2$ is ``almost necessary." Even so, the large deviation principle may still hold with a different rate function and/or a different speed as $p\gg n$ whereas the condition ``$p\gg n^2$" does not hold. We leave it as a future work.

%\section{Proofs}
\section{Proofs of Results in Section \ref{convergence}}\lbl{Proofconvergence}
We start with a concentration inequality on the $\beta$-Hermite ensembles.
\begin{lemma}\lbl{tennis} Suppose $\lambda=(\lambda_1, \cdots, \lambda_n)$ has the joint density function $f_{\beta}(\lambda)$ as in (\ref{betaHermite}). Then, there exists a constant $C>0$ depending on $\beta$ only such that
\beaa
P\big(\max_{1\leq i \leq n}|\lambda_i|\geq \sqrt{n}\,t\big) \leq C\cdot e^{-\frac{1}{2}nt^2 + C n t}
\eeaa
for all $t>0$, $n\geq 2$ and $\beta>0.$
\end{lemma}
Ben Arous et al (2001) study the above probability for the case $\beta=1$ in their Lemma 6.3. Our Lemma \ref{tennis} is stronger than theirs when $t$ is large. In fact, their bound of the above probability is $e^{-\delta nt^2}$ with some $\delta \in (0, 1/2).$\\

\noindent\textbf{Proof}. It is easy to see that the order statistic $\lambda_{(1)}>\cdots > \lambda_{(n)}$ has density function $h_{\beta}(\lambda_1, \cdots, \lambda_n):=n!f_{\beta}(\lambda_1, \cdots, \lambda_n)$ for $\lambda_1>\cdots > \lambda_n.$ Further, $\max_{1\leq i \leq n}|\lambda_i|=|\lambda_{(1)}| \vee |\lambda_{(n)}|.$ It follows that
\beaa
& & P\big(\max_{1\leq i \leq n}|\lambda_i|\geq \sqrt{n}\,t\big)\\
&= & n!K_n^{\beta} \cdot\int_{\lambda_1>\cdots >\lambda_n;\, |\lambda_1| \vee |\lambda_n|\geq \sqrt{n}\,t}\prod_{1\leq i <  j \leq n}|\lambda_i-\lambda_j|^{\beta}\cdot e^{-\frac{1}{2}\sum_{i=1}^n\lambda_i^2}\,d\lambda_1\cdots d\lambda_n\\
& = & \frac{n!K_n^{\beta}}{(n-2)!K_{n-2}^{\beta}}\cdot\int_{\lambda_1>\cdots >\lambda_n;\, |\lambda_1| \vee |\lambda_n|\geq \sqrt{n}\,t}\prod_{i=1}^{n-1}|\lambda_i-\lambda_n|^{\beta}\cdot \prod_{j=2}^{n-1}|\lambda_j-\lambda_1|^{\beta}\cdot e^{-\frac{1}{2}(\lambda_1^2+\lambda_n^2)}\,d\lambda_1d\lambda_n\cdot\\
& &\ \ \ \ \ \ \ \ \ \ \ \ \ \ \ \ \ \ \ \ \ \ \ \ \ \ \ \ \ \ \ \ \ \ \ \ \ \ \ \ \ \ \ \ \ \ \ \ \ \ \ \ \ \ \ \ \ \ \ \ \ \ \ \ \ \ \ \ g(\lambda_2, \cdots, \lambda_{n-1})\,d\lambda_2\cdots d\lambda_{n-1}
\eeaa
where $K_n^{\beta}$ is as in (\ref{Hermitecons}) and
\bea\lbl{histogram}
g(\lambda_2, \cdots, \lambda_{n-1})=(n-2)!K_{n-2}^{\beta}\prod_{2\leq i < j\leq n-1}|\lambda_i-\lambda_j|^{\beta}\cdot e^{-\frac{1}{2}\sum_{i=2}^{n-1}\lambda_j^2}
\eea
for $\lambda_2>\cdots >\lambda_{n-1}.$  Notice
\beaa
\prod_{i=1}^{n-1}|\lambda_i-\lambda_n|^{\beta}\cdot \prod_{j=2}^{n-1}|\lambda_j-\lambda_1|^{\beta}
& \leq & (|\lambda_1| + |\lambda_{n}|)^{(2n-3)\beta} \\
& \leq & \big(2(\lambda_1^2 + \lambda_n^2)\big)^{(2n-3)\beta/2} \leq 2^{\beta n}(\lambda_1^2 + \lambda_{n}^2)^{(2n-3)\beta/2}
\eeaa
for $\lambda_1>\cdots > \lambda_n.$ Further, if $|\lambda_1| \vee |\lambda_n|\geq \sqrt{n}\,t$, then $\lambda_1^2 + \lambda_n^2\geq nt^2.$ Therefore,
\bea
& &  P\big(\max_{1\leq i \leq n}|\lambda_i|\geq \sqrt{n}\,t\big) \\
& \leq & n^22^{\beta n}\frac{K_n^{\beta}}{K_{n-2}^{\beta}} \cdot \int_{x^2 + y^2\geq nt^2}(x^2 + y^2)^{\frac{1}{2}(2n-3)\beta} e^{-\frac{1}{2}(x^2+y^2)}\,dx dy\nonumber\\
& & \ \ \ \ \ \ \ \ \ \ \ \ \ \ \ \ \ \ \ \ \ \ \ \ \ \ \ \ \ \ \ \ \ \ \ \ \ \ \ \ \ \ \ \ \ \ \ \ \ \ \ \ \cdot \int_{\lambda_2>\cdots > \lambda_{n-1}}g(\lambda_2, \cdots, \lambda_{n-1})\,d\lambda_2\cdots d\lambda_{n-1}\nonumber\\
& = & n^22^{\beta n}\frac{K_n^{\beta}}{K_{n-2}^{\beta}} \cdot \int_{x^2 + y^2\geq nt^2}(x^2 + y^2)^{\frac{1}{2}(2n-3)\beta} e^{-\frac{1}{2}(x^2+y^2)}\,dx dy\lbl{rich}
\eea
since $g(\lambda_2, \cdots, \lambda_{n-1})$ is a probability density function. By making transform $x=r\cos\theta$ and $y=r\sin \theta$ with $r\geq \sqrt{n}t$ and $\theta\in [0, 2\pi],$ the last integral is equal to
\bea\lbl{goose}
2\pi\int_{\sqrt{n}t}^{\infty}r^{(2n-3)\beta+1}e^{-r^2/2}\,dr=\pi\int_{nt^2}^{\infty}y^{\frac{\beta}{2}(2n-3)}e^{-y/2}\,dy
\eea
by making another transform $y=r^2.$ To compute the last integral, let's consider $I=\int_b^{\infty}y^{\alpha}e^{-y/2}\,dy$ for $b>0$ and $\alpha>0.$ Use $e^{-y/2}=-2(e^{-y/2})'$ and the integration by parts to have
\beaa
b I \leq  \int_b^{\infty}y^{\alpha+1}e^{-y/2}\,dy=2b^{\alpha+1}e^{-b/2} + 2(\alpha+1) I,
\eeaa
which implies that
\bea\lbl{seminar}
I =\int_b^{\infty}y^{\alpha}e^{-y/2}\,dy\leq \frac{2}{b-2\alpha-2}b^{\alpha+1}e^{-b/2}
\eea
if $\alpha>0$ and $b>2\alpha+2.$ Now, suppose $t>\sqrt{4\beta+4},$ then $nt^2 -(\beta(2n-3) + 2)>\frac{1}{2}nt^2.$ By (\ref{seminar}), the right hand side of (\ref{goose}) is bounded by
\bea\lbl{flawless}
\frac{2\pi}{nt^2-\beta(2n-3)-2}(nt^2)^{\frac{\beta}{2}(2n-3)+1}e^{-nt^2/2} \leq (4\pi) (nt^2)^{\frac{\beta}{2}(2n-3)} e^{-nt^2/2}.
\eea
Now we estimate the term $K_n^{\beta}/K_{n-2}^{\beta}$ appeared in (\ref{rich}). By the Stirling formula (see, e.g., p.204 from Ahlfors (1979) or p.368 from Gamelin (2001)),
\begin{eqnarray}\lbl{implication1}
\log\Gamma(z)=z\log z - z -\frac{1}{2}\log z+ \log \sqrt{2\pi}
+\frac{1}{12z} +O\left(\frac{1}{x^3}\right)
\end{eqnarray}
as $x=\mbox{Re}\,(z)\to +\infty.$ It is easy to check that there exists an absolute constant $C>0$ such that $\Gamma(1+x)\geq x^xe^{-Cx}$ for all $x>0.$ Now, recalling $K_n^{\beta}$ as in (\ref{Hermitecons}), we have
\beaa
\frac{K_n^{\beta}}{K_{n-2}^{\beta}}
=  \frac{1}{2\pi}\cdot \frac{\Gamma(1+\frac{\beta}{2})^2}{\Gamma(1+\frac{\beta}{2}(n-1))\cdot\Gamma(1+\frac{\beta}{2}n)}
 \leq  C_{\beta}\Big(\frac{\beta}{2}(n-1)\Big)^{-\beta(n-1)}e^{\beta nC}.
%\leq C_{\beta}'p^{-\beta p} e^{C_{\beta}'p}
\eeaa
Use $n-1\geq n/2$ to have
\beaa
\Big(\frac{\beta}{2}(n-1)\Big)^{-\beta(n-1)} \leq \Big[\big(\frac{2}{\beta}\big)^{\beta} +1\Big]^n\big(\frac{n}{2}\big)^{-\beta(n-1)} \leq e^{C_{\beta}'n}\cdot 4^{\beta n} n^{-\beta n}
\eeaa
since $\big(\frac{n}{2}\big)^{-\beta(n-1)}=(2^{n-1}n)^{\beta}n^{-\beta n} \leq 4^{\beta n}n^{-\beta n}.$ Combining the last two inequalities we see that $K_n^{\beta}/K_{n-2}^{\beta} \leq C_{\beta}''n^{-\beta n} e^{C_{\beta}''n}.$ This  together with (\ref{rich}) and (\ref{flawless}) concludes that, for some constant $\gamma$ depending on $\beta$ only,
\bea
 P\big(\max_{1\leq i \leq n}|\lambda_i|\geq \sqrt{n}\,t\big) & \leq & \gamma n^{\gamma}t^{n\gamma}\cdot\exp\big(-\frac{1}{2}nt^2 + \gamma n\big)\nonumber\\
 & \leq & \gamma\cdot \exp\big(-\frac{1}{2}nt^2 + \gamma n(t+2)\big)\nonumber\\
 & \leq & (2\gamma)\cdot \exp\big(-\frac{1}{2}nt^2 + (2\gamma) nt\big) \lbl{minister}
\eea
for all $t>\sqrt{4\beta+4},$ where the inequality $n^{\gamma}\vee t^{n\gamma}\leq e^{\gamma nt}$ is used in the second inequality. Note that the last term in (\ref{minister}) is increasing in $\gamma>0.$ Set $\gamma'= \gamma+ \sqrt{\beta +1},$ then $2\gamma'n t \geq nt^2$ for all $0\leq t \leq \sqrt{4\beta +4}.$ It follows that
\beaa
\inf_{0\leq t \leq \sqrt{4\beta +4}}\Big\{(2\gamma')\cdot \exp\big(-\frac{1}{2}nt^2 + (2\gamma') nt\big)\Big\} \geq \inf_{0\leq t \leq \sqrt{4\beta +4}}\{2\gamma' \cdot e^{nt^2/2}\} \geq 1.
\eeaa
Therefore, (\ref{minister}) holds with $C=2\gamma'$ for all $t\geq 0$, $\beta>0$ and $n\geq 2.$\ \ \ \ \ \ \ \ $\blacksquare$\\

\noindent\textbf{Proof of Theorem \ref{main}}. Let $f_{n,\beta}(\lambda)=f_{n,\beta}(\lambda_1, \cdots, \lambda_n)$ be as in (\ref{bWishart}). Recall $x_i=\sqrt{\frac{p}{2\beta }}\big(\frac{\lambda_i}{p}-\beta\big)$ for $1\leq i \leq n$ and $\lambda_i=p\beta +\sqrt{2\beta p}x_i.$ By (\ref{bWishart}),  $x=(x_1, \cdots, x_n)$ has density function
\beaa
\tilde{f}_{n,\beta}(x)=\tilde{c}_{n}^{\beta, p}\prod_{1\leq i<j\leq
n}|x_i-x_j|^{\beta}\cdot\prod_{i=1}^n\Big(1+\sqrt{\frac{2}{\beta p}}x_i\Big)^{\frac{\beta}{2}(p-n+1)-1}\cdot
e^{-\sqrt{\frac{\beta p}{2}}\sum_{i=1}^nx_i}
\eeaa
for all $x_i\geq -\sqrt{\frac{\beta p}{2}}$ with $i=1,\cdots, n$; $\tilde{f}_{n,\beta}(x)=0$ otherwise, where
\bea\lbl{tea}
\tilde{c}_{n}^{\beta, p}=c_{n}^{\beta, p}\cdot (2\beta p)^{\frac{1}{4}n(n-1)\beta +\frac{1}{2}n}\cdot e^{-\frac{1}{2}np\beta}\cdot (p\beta)^{\frac{1}{2}n(p-n+1)\beta - n}.
\eea
Let $\mu_1, \cdots, \mu_n$ have density function $f_{\beta}(\mu)$ as in (\ref{betaHermite}). Then, by (\ref{variation}),
\bea
\|\ml{L}(x_1, \cdots, x_n)-\ml{L}(\mu_1, \cdots, \mu_n)\|
&= & \int_{\mathbb{R}^n}|\tilde{f}_{n,\beta}(x)-f_{\beta}(x)|\,dx_1 \cdots dx_n \nonumber\\
& = &  E\Big|\frac{\tilde{f}_{n,\beta}(X)}{f_{\beta}(X)}-1\Big| \lbl{Malone}
\eea
where the random vector $X=(x_1, \cdots, x_n)\in \mathbb{R}^n$ has density function $f_{\beta}(x)$ (replacing $\lambda$ and $\lambda_i$ in (\ref{betaHermite}) by $x$ and $x_i$ accordingly). Now,
\bea
\frac{\tilde{f}_{n,\beta}(X)}{f_{\beta}(X)}=\frac{\tilde{c}_{n}^{\beta, p}}{K_n^{\beta}}\cdot \prod_{i=1}^n\Big(1+\sqrt{\frac{2}{\beta p}}x_i\Big)^{\frac{\beta}{2}(p-n+1)-1}\cdot
e^{-\sqrt{\frac{\beta p}{2}}\sum_{i=1}^nx_i+\frac{1}{2}\sum_{i=1}^nx_i^2} \lbl{Harper}
\eea
for all $x_i\geq -\sqrt{\frac{\beta p}{2}}$ with $i=1,\cdots, n,$ and it is equal to $0,$ otherwise.

Recall the two conditions: (i) $n\to\infty$ and $p=p_n \gg n^3$ and (ii) $n$ is fixed and $p\to\infty.$ In case (i) we choose constant $t_n>0$ for all $n\geq 1$ such that
\bea\lbl{disjoint}
t_n\to \infty\ \ \ \mbox{and}\ \ \ t_n^4\cdot \frac{n^3}{p} \to 0
\eea
as $n\to\infty.$ In case (ii) we choose $t_p$ for all $p\geq 1$ satisfying
\bea\lbl{doing}
t_p\to \infty\ \ \ \mbox{and}\ \ \ t_p^4\cdot \frac{n^3}{p} \to 0
\eea
as $p\to\infty.$ From now on we will only prove the theorem for case (i). The proof for case (ii)  will be carried through  by replacing ``$t_n$" in (\ref{disjoint}) with ``$t_p$" in (\ref{doing}) and ``$n\to\infty$" with ``$p\to\infty$" in the context.

By Lemma \ref{tennis}, $P\big(\max_{1\leq i \leq n}|x_i|\geq \sqrt{n}\,t_n\big) \leq C\cdot e^{-\frac{1}{3}nt_n^2} \to 0$ as $n\to\infty$. Set
\bea\lbl{device}
\Omega_n=\Big\{\max_{1\leq i \leq n}\frac{|x_i|}{\sqrt{n}}\leq t_n\Big\},\ \ n\geq 1.
\eea
Then $P(\Omega_n)\to 1$ as $n\to\infty.$ By the Taylor expansion, there exists $\epsilon_0 \in (0,1)$ such that
\bea\lbl{open}
\log (1+x)= x-\frac{1}{2}x^2 +\frac{1}{3}x^3+ u(x)\ \ \ \mbox{with}\ \ \ \ |u(x)| \leq |x|^4
\eea
for all $|x|<\epsilon_0.$ Notice $\sqrt{\frac{2}{\beta p}}|x_i| <\epsilon_0$ on $\Omega_n$ as $n$ is large enough by (\ref{disjoint}). It follows that
\bea
& & \log \prod_{i=1}^n\Big(1+\sqrt{\frac{2}{\beta p}}x_i\Big)^{\frac{\beta}{2}(p-n+1)-1}\nonumber\\
& = & \Big(\frac{\beta}{2}(p-n+1)-1\Big)\sum_{i=1}^n\log \Big(1+\sqrt{\frac{2}{\beta p}}x_i\Big)\nonumber\\
& = & \Big(\frac{\beta}{2}(p-n+1)-1\Big)\Big(\sqrt{\frac{2}{\beta p}}\,\sum_{i=1}^nx_i-\frac{1}{\beta p}\sum_{i=1}^nx_i^2 + \frac{1}{3}\sqrt{\frac{8}{\beta^3p^3}}\sum_{i=1}^nx_i^3 +\sum_{i=1}^n u\Big(\sqrt{\frac{2}{\beta p}}x_i\Big)\Big)\nonumber\\
& & \ \ \ \ \ \ \ \ \ \ \ \ \ \ \ \ \ \ \ \ \ \ \ \ \ \ \ \ \ \ \ \ \ \ \ \ \  \ \ \ \lbl{shaq}
\eea
on $\Omega_n$ as $n$ is sufficiently large. By writing $\frac{\beta}{2}(p-n+1)-1=\frac{\beta}{2}p-(\frac{\beta(n-1)}{2}+1),$ we have
\bea
U_{n,1}:&=&\Big(\frac{\beta}{2}(p-n+1)-1\Big)\cdot\sqrt{\frac{2}{\beta p}}\,\sum_{i=1}^nx_i\nonumber\\
&=& \sqrt{\frac{\beta p}{2}}\sum_{i=1}^nx_i + \delta_{n,1}\cdot \Big(\frac{n^3}{p}\Big)^{1/2}\Big(\sum_{i=1}^n\frac{x_i}{\sqrt{n}}\Big) \lbl{kobe}
\eea
with $|\delta_{n,1}| \leq C_{\beta,1},$ where $C_{\beta,1}>0$ is a constant depending on $\beta$ only. Second,
\bea
U_{n,2}:=-\Big(\frac{\beta}{2}(p-n+1)-1\Big)\cdot \frac{1}{\beta p}\sum_{i=1}^nx_i^2
&=& -\frac{1}{2}\sum_{i=1}^nx_i^2+ \delta_{n,2}\cdot \Big(\frac{ n}{p}\sum_{i=1}^nx_i^2\Big)\nonumber\\
& = & -\frac{1}{2}\sum_{i=1}^nx_i^2+ \delta_{n,2}\cdot \frac{ n^3t_n^2}{p}\lbl{Jordan}
\eea
on $\Omega_n$ with $|\delta_{n,2}| \leq C_{\beta,2},$ where $C_{\beta,2}>0$ is a constant depending on $\beta$ only. Now,
\bea
U_{n,3}:& = &\Big(\frac{\beta}{2}(p-n+1)-1\Big)\cdot\frac{1}{3}\sqrt{\frac{8}{\beta^3p^3}}\sum_{i=1}^nx_i^3 \nonumber\\
& = & \delta_{n,3}\cdot\Big(\frac{n^3}{p}\Big)^{1/2}\sum_{i=1}^n\Big(\frac{x_i}{\sqrt{n}}\Big)^3 \lbl{Pippen}
\eea
with $|\delta_{n,3}| \leq C_{\beta,3},$ where $C_{\beta,3}>0$ is a constant depending on $\beta$ only. On $\Omega_n,$ by (\ref{open})
\bea
U_{n,4}:=  \Big(\frac{\beta}{2}(p-n+1)-1\Big)\cdot \Big|\sum_{i=1}^n u\Big(\sqrt{\frac{2}{\beta p}}x_i\Big)\Big|
& \leq & \delta_{n,4}\cdot\Big(\frac{1}{p}\sum_{i=1}^n|x_i|^4\Big) \nonumber\\
& \leq & \delta_{n,4}\cdot \frac{n^3t_n^4}{p}\lbl{Rodman}
\eea
with $|\delta_{n,4}| \leq C_{\beta,4},$ where $C_{\beta,4}>0$ is a constant depending on $\beta$ only. We claim that
\bea\lbl{cumulant}
& & U_{n,1}=\sqrt{\frac{\beta p}{2}}\sum_{i=1}^nx_i +o_P(1),\ \ U_{n,2}= -\frac{1}{2}\sum_{i=1}^nx_i^2+ o_P(1), \nonumber\\
& & U_{n,3}=o_P(1), \ \ \ U_{n,4}=o_P(1)
\eea
as $n\to\infty,$ where by $Z_n=o_P(1)$ we mean $Z_n\to 0$ in probability as $n\to\infty.$

Looking at (\ref{Jordan}) and (\ref{Rodman}) together with (\ref{disjoint}) and the fact $P(\Omega_n)\to 1$, the claims for $U_{n,2}$ and $U_{n,4}$ in (\ref{cumulant}) are obviously true. By Theorem 1.2 from Dumitriu and Edelman (2006),  for each integer $k\geq1$,  $\sum_{i=1}^n(x_i/\sqrt{n})^{2k-1}$ converges in distribution to $N(0, \sigma_k^2)$ with $\sigma_k^2<\infty$ as $n\to\infty$. Reviewing (\ref{kobe}) and (\ref{Pippen}), from the condition $p \gg n^3$, (\ref{disjoint}) and the fact $P(\Omega_n)\to 1,$ we see that the claims for $U_{n,2}$ and $U_{n,4}$ in (\ref{cumulant}) hold true. (If $p_n \equiv p\geq 2$ for all $n\geq 1$ then the claims for  $U_{n,1}$ and $U_{n,3}$ are evidently true by (\ref{kobe}) and (\ref{Pippen}) together with the fact $P(\Omega_n)\to 1$).

Now, from (\ref{shaq})-(\ref{Rodman}), we have
\beaa
\log \prod_{i=1}^n\Big(1+\sqrt{\frac{2}{\beta p}}x_i\Big)^{\frac{\beta}{2}(p-n+1)-1}=U_{n,1} + U_{n,2} +U_{n,3} +U_{n,4}
\eeaa
on $\Omega_n.$ Consequently,  from (\ref{Harper}), (\ref{cumulant}) and the fact $P(\Omega_n)\to 1$, we conclude that
\bea\lbl{Gasol}
\frac{\tilde{f}_{n,\beta}(X)}{f_{\beta}(X)}\cdot\Big(\frac{\tilde{c}_{n}^{\beta, p}}{K_n^{\beta}}\Big)^{-1}\to 1
\eea
in probability as $n\to\infty.$ Next we show
\bea\lbl{google}
\lim_{n\to\infty}\frac{\tilde{c}_{n}^{\beta, p}}{K_n^{\beta}} = 1.
\eea

%\begin{lemma}\lbl{river} Given $\beta>0$. Let $\tilde{c}_{n}^{\beta, p}$ be as in (\ref{tea}) and $K_p^{\beta}$ be as in %(\ref{Hermitecons}). If $p=o(n^{1/3}),$ then $ as $n\to\infty.$
%\end{lemma}
%\noindent\textbf{Proof of Lemma \ref{river}}.

\noindent Recall (\ref{conWishart}) and (\ref{Hermitecons}), we know
\bea\lbl{divider}
\frac{c_{n}^{\beta,p}}{K_n^{\beta}} = (2\pi)^{n/2}\cdot 2^{-\beta np/2}\prod_{j=1}^{n}\frac{1}
{\Gamma(\frac{\beta}{2}(p-n+j))}= (2\pi)^{n/2}\cdot 2^{-\beta np/2}\prod_{i=0}^{n-1}\frac{1}
{\Gamma(\frac{\beta}{2}(p-i))}.
\eea
From (\ref{implication1}) we have
\bea
& & \log \Gamma(\frac{\beta}{2}(p-i))\nonumber\\
&= &\frac{\beta}{2}(p-i)\log \Big(\frac{\beta}{2}(p-i)\Big) - \frac{\beta}{2}(p-i) -\frac{1}{2}\log \big(\frac{\beta}{2}(p-i)\big) +\log \sqrt{2\pi} + O(\frac{1}{p-i})\nonumber\\
& = & \frac{\beta}{2}(p-i)\log (p-i) + (\frac{\beta}{2}\log \frac{\beta}{2}-\frac{\beta}{2})(p-i)-\frac{1}{2}\log (p-i) -\frac{1}{2}\log \frac{\beta}{2}+\log \sqrt{2\pi}+ O(\frac{1}{p})\nonumber\\
& & \ \ \ \ \ \ \ \ \ \ \ \ \ \ \ \ \ \ \ \ \ \ \ \ \ \  \lbl{whenever}
\eea
as $n\to\infty$ uniformly for all $i=0,1,\cdots, n-1.$ Use the condition $p \gg n$ and the fact
\bea\lbl{cheese}
\log (p-i)=\log p -\frac{i}{p} +O\Big(\frac{n^2}{p^2}\Big)
\eea
uniformly for all $0\leq i \leq n-1$ as $n\to \infty$ to have
\bea\lbl{stunning}
-\frac{1}{2}\sum_{i=0}^{n-1}\log (p-i)=-\frac{1}{2}n\log p + O(\frac{n^2}{p})
\eea
as $n\to\infty.$ Moreover,
\bea\lbl{factoring}
(\frac{\beta}{2}\log \frac{\beta}{2}-\frac{\beta}{2})\sum_{i=0}^{n-1} (p-i) = \frac{1}{2}(\frac{\beta}{2}\log \frac{\beta}{2}-\frac{\beta}{2})n(2p-n+1).
%\nonumber\\
%& = & \frac{1}{2}(\frac{\beta}{2}\log \frac{\beta}{2}-\frac{\beta}{2})p(2n-p) + O(p)
\eea
 By (\ref{cheese}) again,
\beaa
\frac{\beta}{2}\sum_{i=0}^{n-1}(p-i)\log (p-i)
& = & \Big[\frac{\beta}{2}\sum_{i=0}^{n-1}(p-i)(\log p -\frac{i}{p})\Big] + O(\frac{n^3}{p})\\
& = & \frac{\beta}{4}(\log p)n(2p-n+1)- \Big[\frac{\beta}{2}\sum_{i=0}^{n-1}(i-\frac{i^2}{p})\Big]+O(\frac{n^3}{p})\\
& = & \frac{\beta}{4}(\log p)n(2p-n+1)-\frac{\beta}{4}n(n-1) +O(\frac{n^3}{p})
\eeaa
as $n\to\infty,$ where we use the fact $\sum_{i=0}^{n-1}(i-\frac{i^2}{p})=\frac{1}{2}n(n-1)+O(\frac{n^3}{p})$ in the last step. This joint with (\ref{whenever}), (\ref{stunning}) and (\ref{factoring}) leads to
\bea\lbl{lemonade}
& & \log \Big(\prod_{j=1}^n\Gamma(\frac{\beta}{2}(p-n+j))\Big)^{-1} \nonumber\\
& = & -\sum_{i=0}^{n-1}\log \Gamma(\frac{\beta}{2}(p-i)) \nonumber\\
& = & -\frac{\beta}{4}(\log p)n(2p-n+1)- \frac{\beta}{2}(\log \frac{\beta}{2}-1)np + \frac{\beta}{4}(\log \frac{\beta}{2})n^2+\frac{1}{2}n\log p + O\big(n+\frac{n^3}{p}\big)\nonumber\\
& &
\eea
as $n\to\infty$ under the restriction that $p\gg n$ only. Consequently,
\bea
\log \frac{c_{n}^{\beta,p}}{K_n^{\beta}}&= & \frac{n}{2}\log (2\pi)-\frac{\beta np}{2}\log 2-\sum_{i=0}^{n-1}\log \Gamma(\frac{\beta}{2}(p-i)) \nonumber\\
&= & -\frac{\beta np}{2}\log 2 +\big(\frac{1}{2}\log \frac{\beta}{2}\big)n +\frac{1}{2}n\log p - \frac{1}{2}(\frac{\beta}{2}\log \frac{\beta}{2}-\frac{\beta}{2})n(2p-n+1) \nonumber\\
& & \ \ \ \ \ \ \ \ \ \ \ \ \ \ \   - \frac{\beta}{4}(\log p)n(2p-n+1)+\frac{\beta}{4}n(n-1) +O(\frac{n^3}{p}) \lbl{gull}
\eea
as $n\to\infty.$ Reviewing (\ref{tea}),
\bea\lbl{sound}
& & \log \frac{\tilde{c}_{n}^{\beta,p}}{K_n^{\beta}}=\log \frac{c_{n}^{\beta,p}}{K_n^{\beta}}+ \big(\frac{1}{4}n(n-1)\beta +\frac{1}{2}n\big)\log (2\beta p)-\frac{1}{2}np\beta \nonumber\\
& & \ \ \ \ \ \ \ \ \ \ \ \ \ \ \ \ \ \ \ \ \ \ \ \  + \big(\frac{1}{2}n(p-n+1)\beta - n\big)\log (p\beta ).
\eea
Combining this with (\ref{gull}), by a routine but tedious calculation (see it in Appendix), we have
\bea\lbl{oversize}
\log \frac{\tilde{c}_{n}^{\beta,p}}{K_n^{\beta}}=O(\frac{n^3}{p})\to 0
\eea
as $n\to\infty,$ which implies (\ref{google}). Finally, by (\ref{Gasol}) and (\ref{google}),
\beaa
\frac{\tilde{f}_{n,\beta}(X)}{f_{\beta}(X)}\to 1
\eeaa
in probability as $n\to\infty.$ Obviously, $E\frac{\tilde{f}_{n,\beta}(X)}{f_{\beta}(X)}=1$ for all $n\geq 2.$ By a variant of the Scheff\'{e} Lemma (see, e.g., Corollary 4.2.4 from Chow and Teicher (1997)), the two facts imply that $E|\frac{\tilde{f}_{n,\beta}(X)}{f_{\beta}(X)}-1|\to 0$ as $n\to\infty.$ The desired conclusion then follows from (\ref{Malone}). \ \ \ \ \ \ \ \ $\blacksquare$\\

\noindent\textbf{Proof of Proposition \ref{independence}}. Let $\xi_1, \cdots, \xi_n$ have density function $f_{2}(\xi_1, \cdots, \xi_n)$ as in (\ref{betaHermite}) with $\beta=2$. Then $y_1:=\xi_1/\sqrt{2}, \cdots, y_n:=\xi_n/\sqrt{2}$ have density function
\beaa
f(y_1, \cdots, y_n)= \mbox{Const}\cdot \prod_{1\leq i<j\leq
n}|y_i-y_j|^{2}\cdot
e^{-\sum_{i=1}^ny_i^2}
\eeaa
for $(y_1, \cdots, y_n)\in \mathbb{R}^n.$ It is shown by Bornemann (2010) (see also Bianchi et al (2010))  that the two random variables $\tilde{y}_{min}:=\sqrt{2}n^{1/6}(y_{min} +\sqrt{2n})=n^{1/6}(\xi_{min} +2\sqrt{n})$ and $\tilde{y}_{max}:=\sqrt{2}n^{1/6}(y_{max} -\sqrt{2n})=n^{1/6}(\xi_{max} -2\sqrt{n})$ are asymptotic independent, that is,
\bea\lbl{mother}
P\big( \tilde{y}_{min}\in A,\, \tilde{y}_{max} \in B\big) - P\big(\tilde{y}_{min} \in A\big)\cdot P\big(\tilde{y}_{max} \in B\big) \to 0
\eea
as $n\to\infty$ for any Borel sets $A$ and $B.$  Further, $\tilde{y}_{max}$ goes weakly to $U$ and $\tilde{y}_{min}$ goes weakly to $-V,$ where $U$ are $V$ are i.i.d. with the distribution function $F_2(x)$ as in (\ref{F2}) (see also Tracy and Widom (1993, 1994)).  By the assumptions, $\lambda=(\lambda_1, \cdots, \lambda_{n})$ has density function $f_{n,2}(\lambda)$ as in (\ref{bWishart}). In (\ref{mother}) replacing $\xi_{min}$ in the expression of $\tilde{y}_{min}$ by  $\frac{\sqrt{p}}{2}\big(\frac{\lambda_{min}}{p}-2\big)$ and $\xi_{max}$ in the expression of $\tilde{y}_{max}$  by $\frac{\sqrt{p}}{2}\big(\frac{\lambda_{max}}{p}-2\big)$, respectively, we obtain from  Theorem \ref{main} that
\bea
& & P\Big( \frac{\lambda_{min}-\mu_{n,1}}{\sigma_n}\in A,\, \frac{\lambda_{max}-\mu_{n,2}}{\sigma_n} \in B\Big) - P\Big( \frac{\lambda_{min}-\mu_{n,1}}{\sigma_n}\in A\Big) \cdot \Big(\frac{\lambda_{max}-\mu_{n,2}}{\sigma_n} \in B\Big) \ \ \  \nonumber\\
& & \to 0\lbl{cough}
\eea
as $n\to\infty,$ where $\mu_{n,1}=2p-4\sqrt{pn}$, $\mu_{n,2}=2p+4\sqrt{pn}$ and $\sigma_n=2\sqrt{p}n^{-1/6}.$ That is, $(\lambda_{min}-\mu_{n,1})/\sigma_n$ and $(\lambda_{max}-\mu_{n,2})/\sigma_n$ are asymptotic independent.

Finally, using the same argument as in the above, the weak convergence of $\tilde{y}_{max}$ to $U$ and that of $\tilde{y}_{min}$ to $-V,$ we obtain that $\frac{\lambda_{max}-\mu_{n,2}}{\sigma_n}$ converges weakly to $U$ and $\frac{\lambda_{min}-\mu_{n,1}}{\sigma_n}$ converges weakly to $-V$ as $n\to\infty.$ This together with (\ref{cough}) gives the desired conclusion.\ \ \ \ \ \ \ $\blacksquare$\\

\noindent\textbf{Proof of Corollary \ref{prasing}.}  Let $\lambda_1, \cdots, \lambda_{n}$ be the eigenvalues of $\bd{X}\bd{X}^*.$ As mentioned before (\ref{bWishart}), we know $\lambda=(\lambda_1, \cdots, \lambda_{n})$ has density function $f_{n,\beta}(\lambda)$ as in (\ref{bWishart}) with $\beta=2$. Recall $\mu_{n,1}=2p-4\sqrt{pn}$, $\mu_{n,2}=2p+4\sqrt{pn}$ and $\sigma_n=2\sqrt{p}n^{-1/6}$ in Proposition \ref{independence}. Since $\sigma_n\to \infty$ and $\mu_{n,1}/(2p) \to 1$, by the Slusky lemma, $\lambda_{min}/(2p)\to 1$ in probability as $n\to\infty.$ Set $\delta_n=4\sqrt{\frac{n}{p}}.$ Write
\beaa
& & \sqrt{p}n^{1/6}\Big(\frac{\lambda_{max}}{\lambda_{min}}-1-4\sqrt{\frac{n}{p}}\Big)\\
& = & \frac{2p}{\sigma_n}\cdot \frac{\lambda_{max}-(1+\delta_n)\lambda_{min}}{\lambda_{min}}\\
& = & \frac{2p}{\lambda_{min}}\cdot\Big(\frac{\lambda_{max}-\mu_{n,2}}{\sigma_n}- (1+\delta_n)\frac{\lambda_{min}-\mu_{n,1}}{\sigma_n}+ \frac{\mu_{n,2}-(1+\delta_n)\mu_{n,1}}{\sigma_n}\Big).
\eeaa
It is easy to check the last term in the parenthesis is equal to $8n^{7/6}p^{-1/2}\to 0$ since $p \gg n^3.$ Also, $1+\delta_n\to 1$ and $\frac{2p}{\lambda_{min}}\to 1$ in probability as $n\to\infty.$ These and Proposition \ref{independence} conclude that
\bea\lbl{up}
\sqrt{p}n^{1/6}\Big(\frac{\lambda_{max}}{\lambda_{min}}-1-4\sqrt{\frac{n}{p}}\Big)
\eea
converges weakly to $U+V,$ where $U$ and $V$ are i.i.d. with distribution function $F_2(x)$ as in (\ref{F2}). Now, let  $\alpha_n=2\sqrt{p}n^{1/6}$, $\beta_n=1+2\sqrt{\frac{n}{p}}.$ Recall $\kappa_n$ defined in (\ref{Changchun}). Observe
\beaa
\alpha_n(\kappa_n -\beta_n) & = & \alpha_n\Big(\frac{\lambda_{max}}{\lambda_{min}} -\beta_n^2\Big)\cdot \big(\sqrt{\frac{\lambda_{max}}{\lambda_{min}}} +\beta_n\Big)^{-1}\\
& = & \Big\{\sqrt{p}n^{1/6}\Big(\frac{\lambda_{max}}{\lambda_{min}}-1-4\sqrt{\frac{n}{p}}\Big)- \frac{4n^{7/6}}{\sqrt{p}}\Big\}\cdot \big(\sqrt{\frac{\lambda_{max}}{\lambda_{min}}} +\beta_n\Big)^{-1}\cdot 2.
\eeaa
From (\ref{up}), we see that $\lambda_{max}/\lambda_{min}\to 1$ in probability as $n\to\infty.$ Also, $\beta_n\to 1$ and  $4n^{7/6}/\sqrt{p}\to 0$  since $p\gg n^3.$ The desired conclusion then follows from  the Slusky lemma and (\ref{up}). \ \ \ \ \ \ \ $\blacksquare$\\

\noindent\textbf{Proof of Proposition \ref{singing}}. Let $\xi_1, \cdots, \xi_n$ have density function $f_{\beta}(\xi_1, \cdots, \xi_n)$ as in (\ref{betaHermite}). Then $y_i:=\sqrt{\frac{2}{\beta}}\,\xi_i,\, 1\leq i \leq n,$ have density function
\beaa
f(y_1, \cdots, y_n)= \mbox{Const}\cdot \prod_{1\leq i<j\leq
n}|y_i-y_j|^{\beta}\cdot
e^{-\frac{\beta}{4}\sum_{i=1}^n y_i^2}
\eeaa
for $(y_1, \cdots, y_n)\in \mathbb{R}^n.$ Since $-\xi_{min}$ and $\xi_{max}$ have the same distribution, by Theorem 1.1 from Ram\'{i}rez et al (2011), $n^{1/6}(2\sqrt{n}+\sqrt{2/\beta}\,\xi_{min}) \overset{d}{=} n^{1/6}(2\sqrt{n}-\sqrt{2/\beta}\,\xi_{max}) $ converges weakly to the distribution of $\Lambda_0,$ which is defined below (\ref{radish}). Let $\lambda=(\lambda_1, \cdots, \lambda_{n})$ have density function $f_{n,\beta}(\lambda)$ as in (\ref{bWishart}). By  Theorem \ref{main},
\beaa
P\Big(g_n\Big(\sqrt{\frac{p}{2\beta }}\big(\frac{\lambda_{min}}{p}-\beta\big)\Big) \leq x\Big)- P\big(g_n(\xi_{min})\leq x\big) \to 0
\eeaa
as $n\to\infty$ for any   $x\in \mathbb{R}$ and any sequence of Borel measurable functions  $\{g_n(t);\, t\in \mathbb{R},\, n\geq 2\}.$  Taking $g_n(t)=n^{1/6}(2\sqrt{n}+\sqrt{2/\beta}\,t)$ to get
\beaa
P\Big(\frac{\lambda_{min}-\mu_n}{\sigma_n}\leq x\Big)- P\Big(n^{1/6}\Big(2\sqrt{n}+\sqrt{\frac{2}{\beta}}\,\xi_{min}\Big)\leq x\Big) \to 0
\eeaa
where $\mu_n=\beta (p-2 \sqrt{np})$ and $\sigma_n=\beta \sqrt{p}n^{-1/6}.$ By the earlier conclusion, the last probability goes to $H(x)=P(\Lambda_0\leq x)$ for all continuous point $x$ of $H(x).$ This leads to the desired conclusion.\ \ \ \ \ \ \ $\blacksquare$

\section{Proofs of Results in Section \ref{LDPsection}}\lbl{Proofrainbow}
This section is divided into two subsections. In Subsection \ref{Proof_pack} we prove Theorems \ref{pack} and \ref{done1} for the large deviations for the largest and the smallest eigenvalues of the $\beta$-Laguerre ensembles.
%Subsection \ref{Proof_pack} is devoted to the proof of Theorem \ref{pack}: the large deviations for the largest eigenvalues of the $\beta$-Laguerre ensembles.
Subsection \ref{ProofLDP} is devoted to the proof of  Theorem \ref{main1} for the large deviations for the empirical distributions of the eigenvalues from the same ensembles.

\subsection{Proof of Theorem \ref{pack}}\lbl{Proof_pack}
%\setcounter{equation}{0}
%\section{Some Auxiliary Results on $\beta$-Laguerre Ensembles}\lbl{sectionLaguerre}
%\setcounter{equation}{0}
%We first estimate the product of a sequence of Gamma functions, which will be used several times later.

\noindent\textbf{Proof of Theorem \ref{pack}}. It is easy to check that  $I(x)> 0$ for all $x\ne \beta$, $I(\beta)=0$, $\{I(x)\leq c\}$ is compact for any $c\geq 0,$ and $I(x)$ is strictly increasing on $[\beta, \infty).$ Now, to prove the theorem, we need to show the following
\bea
& & \limsup_{n\to\infty}\frac{1}{p}\log P\Big(\frac{\lambda_{max}}{p} \in F\Big) \leq -\inf_{x\in F} I(x)\lbl{graffiti}\\
& & \liminf_{n\to\infty}\frac{1}{p}\log P\Big(\frac{\lambda_{max}}{p} \in G\Big) \geq -\inf_{x\in G} I(x)\lbl{deer}
\eea
for any closed set $F\subset \mathbb{R}$ and open set $G \subset \mathbb{R}.$

\noindent{{\it The proof of (\ref{graffiti}).}} Obviously, the joint density function of the order statistics $\lambda_{max}=\lambda_{(1)}>\cdots > \lambda_{(n)}$ is $g(\lambda_1, \cdots, \lambda_n)=n!f_{n,\beta}(\lambda_{1}, \cdots, \lambda_n)$ for all $\lambda_1>\cdots >\lambda_n.$ Write
\bea\lbl{stone}
g(\lambda_1, \cdots, \lambda_n)=\underbrace{\frac{n!c_{n}^{\beta, p}}{(n-1)!c_{n-1}^{\beta, p-1}}}_{A_n}\cdot \Big(\underbrace{\lambda_{1}^{\frac{\beta}{2}(p-n+1)-1}
e^{-\frac{1}{2}\lambda_{1}}\prod_{i=2}^{n}(\lambda_{1}-\lambda_i)^{\beta}}_{B_n}\Big)\cdot L_n(\lambda_2, \cdots, \lambda_n)
\eea
where
\bea\lbl{holy}
L_n(\lambda_2, \cdots, \lambda_n)=(n-1)!c_{n-1}^{\beta, p-1}\prod_{2\leq i<j\leq
n}|\lambda_i-\lambda_j|^{\beta}\cdot\prod_{i=2}^n\lambda_i^{\frac{\beta}{2}(p-n+1)-1}\cdot
e^{-\frac{1}{2}\sum_{i=2}^n\lambda_i}.
\eea
Notice that $\prod_{i=2}^{n}|\lambda_{1}-\lambda_i|^{\beta} \leq \lambda_1^{\beta (n-1)}$ for all $\lambda_1>\cdots >\lambda_n.$ This gives
\beaa
B_n \leq \lambda_1^{\frac{\beta}{2}(n+p-1)-1}e^{-\lambda_1/2}.
\eeaa
Thus, from (\ref{stone}) we have
\bea
P\Big(\frac{\lambda_{max}}{p} \geq x\Big) & = &  \int_{px<\lambda_1,\,  \lambda_1>\cdots >\lambda_n>0}g(\lambda_1, \cdots, \lambda_n)\,d\lambda_1\cdots \,d\lambda_n \nonumber\\
& \leq & A_n \cdot \int_{px}^{\infty}\lambda_1^{\frac{\beta}{2}(n+p-1)-1}e^{-\lambda_1/2}\,d\lambda_1 \cdot \int_{\lambda_2>\cdots >\lambda_n>0}L_n(\lambda_2, \cdots, \lambda_n)\,d\lambda_2\cdots \,d\lambda_n\nonumber\\
& = &  A_n \cdot \int_{px}^{\infty}y^{\frac{\beta}{2}(n+p-1)-1}e^{-y/2}\,dy\lbl{here}
\eea
since $L_n(\lambda_2, \cdots, \lambda_n)$  is a probability density function. We claim, as $n$ is sufficiently large, the following hold:
\bea
& & \int_{px}^{\infty}y^{\frac{\beta}{2}(n+p-1)-1}e^{-y/2}\,dy \leq \frac{2}{p(x-\beta)-\beta n +\beta -2}(p x)^{\frac{\beta}{2}(n+p-1)} e^{-px/2}\ \ \mbox{if $x>\beta$}\lbl{hit};\ \ \ \ \\
& & \int_0^{px}y^{\frac{\beta}{2}(p-n+1)-1}e^{-y/2}\,dy \leq \frac{2}{(\beta -x)p-\beta n}(p x)^{\frac{\beta}{2}(p-n+1)} e^{-px/2}\ \ \mbox{if $0<x<\beta$}.\lbl{use}
\eea
In fact, taking $\alpha=\frac{\beta}{2}(n+p-1)-1$ and $b=px$ in (\ref{seminar}) and using the fact $x>\beta,$ we obtain (\ref{hit}).
%\bea\lbl{hit}
%\int_{nx}^{\infty}y^{\frac{\beta}{2}(n+p-1)-1}e^{-y/2}\,dy \leq \frac{2}{n(x-\beta)-\beta p +\beta -2}(n x)^{\frac{\beta}{2}(n+p-1)} %e^{-nx/2}
%\eea
%as $n$ is sufficiently large.
To prove (\ref{use}), set $J=\int_0^{b}y^{\alpha-1}e^{-y/2}\,dy$ with $\alpha=\frac{\beta}{2}(p-n+1)$, $b=px$ and $0<x<\beta.$ By integration by parts,
\beaa
\alpha J=\int_0^{b}(y^{\alpha})'e^{-y/2}\,dy=b^{\alpha}e^{-b/2}+\frac{1}{2}\int_0^{b}y^{\alpha}e^{-y/2}\,dy \leq b^{\alpha}e^{-b/2}+\frac{b}{2}J.
\eeaa
Solve the inequality to have
\beaa
J \leq \frac{2}{2\alpha -b}b^{\alpha}e^{-b/2}\leq \frac{2}{(\beta -x)p-\beta n}(p x)^{\frac{\beta}{2}(p-n+1)} e^{-px/2},
\eeaa
which leads to (\ref{use}).

Now we estimate $A_n$ in (\ref{stone}). In fact, by (\ref{conWishart}),
\bea\lbl{train}
A_n=n\frac{c_{n}^{\beta, p}}{c_{n-1}^{\beta, p-1}}=\frac{n2^{-\beta (n+p-1)/2}\Gamma(1+\frac{\beta}{2})}{\Gamma(1+\frac{\beta}{2}n)\Gamma(\frac{\beta}{2}p)}.
%\cdot \frac{\prod_{j=1}^{p-1}\Gamma(\frac{\beta}{2}(n-p+j+1))}{\prod_{j=1}^{p-1}\Gamma(\frac{\beta}{2}(n-p+j))}.
\eea
%Observe the last ratio is identical to
%\beaa
%\frac{\prod_{j=2}^{p}\Gamma(\frac{\beta}{2}(n-p+j))}{\prod_{j=1}^{p-1}\Gamma(\frac{\beta}{2}(n-p+j))}
%=\frac{\Gamma(\frac{\beta}{2}n)}{\Gamma(\frac{\beta}{2}(n-p+1))}.
%\eeaa
%Thus,
%\beaa
%A_n=\frac{p2^{-\beta n/2}\Gamma(1+\frac{\beta}{2})}{\Gamma(1+\frac{\beta}{2}p)\Gamma(\frac{\beta}{2}(n-p+1))}.
%\eeaa
Use the fact $\Gamma(x+1)=x\Gamma(x)$ to have
\beaa
A_n=\frac{2}{\beta}\cdot\frac{2^{-\beta (n+p-1)/2}\Gamma(1+\frac{\beta}{2})}{\Gamma(\frac{\beta}{2}n)\Gamma(\frac{\beta}{2}p)}.
\eeaa
%Recall the Stirling formula (see, e.g., p.368 from \cite{Gamelin} or (37) on
%p.204 from \cite{Ahlfors}):
%\begin{eqnarray}\lbl{implication}
%\log\Gamma(z)=z\log z - z -\frac{1}{2}\log z+ \log \sqrt{2\pi}
%+\frac{1}{12z} +O\left(\frac{1}{x^3}\right)
%\end{eqnarray}
%as $x=\mbox{Re}\,(z)\to +\infty$.
By (\ref{implication1}),
\bea
\log A_n &= & -\frac{\beta}{2} (n+p-1)\log 2- \log \Gamma(\frac{\beta}{2}n)-\log \Gamma(\frac{\beta}{2}p)+o(p)\nonumber\\
%&  &  +\frac{\beta}{2}(n-p+1)-\frac{\beta}{2}(n-p+1)\log \big[\frac{\beta}{2}(n-p+1)\big] + o(n)\\
&= &-\frac{\beta}{2}(\log 2) p-\frac{\beta}{2}n\log n +\frac{\beta}{2}p -\frac{\beta}{2}p\log (\frac{\beta}{2}p)+ o(p)\nonumber\\
& = & -\frac{\beta}{2}p\big(\log 2 -1+\log \frac{\beta}{2}\big) -\frac{\beta}{2}n\log n -\frac{\beta}{2}p\log p +o(p)\lbl{outrage}
\eea
as $n\to\infty.$ Combining (\ref{here}) and (\ref{hit}) we have
\bea
\log P\Big(\frac{\lambda_{max}}{p} \geq x\Big) &\leq &  -\frac{\beta}{2}p\big(\log 2 -1+\log \frac{\beta}{2}\big) -\frac{\beta}{2}n\log n -\frac{\beta}{2}p\log p \nonumber\\
 & & -\frac{p x}{2} + \frac{\beta}{2}(p+n-1)\log p + \frac{\beta}{2}(p+n-1)\log x + o(p) \nonumber\\
 & = & p\big(\frac{\beta}{2}-\frac{\beta}{2}\log \beta-\frac{x}{2}+\frac{\beta}{2}\log x\big)
  -\frac{\beta}{2}n\log n +\frac{\beta}{2}n\log p+o(p)\nonumber\\
  & & \lbl{cabin}
\eea
as $n\to\infty$ for all $x>\beta.$ Since $-\frac{\beta}{2}n\log n +\frac{\beta}{2}n\log p=-\frac{\beta}{2}n\log \frac{n}{p}=o(p)$ as $n\to\infty,$  we arrive at
%Note that $(p\log p)\vee(p\log n)=o(n)$ since $p=o(n/\log n).$ This joint with the above concludes that
\bea\lbl{compared1}
\limsup_{n\to\infty}\frac{1}{p}\log P\Big(\frac{\lambda_{max}}{p} \geq x\Big) \leq -\Big(\frac{x-\beta}{2}-\frac{\beta}{2}\log \frac{x}{\beta}\Big)
\eea
for any $x>\beta,$ and hence the same holds for $x\geq \beta$ since the right hand side of (\ref{compared1}) is equal to zero when $x=\beta.$

Now, if $\frac{1}{p}\lambda_{max}\leq x\in (0, \beta),$ by the definition of $f_{n,\beta}(\lambda)$ in (\ref{bWishart}), we see
\beaa
f_{n,\beta}(\lambda)\leq c_{n}^{\beta, p}(p x)^{\beta n(n-1)/2}\cdot \prod_{i=1}^n\lambda_i^{\frac{\beta}{2}(p-n+1)-1}\cdot
e^{-\frac{1}{2}\sum_{i=1}^n\lambda_i}
\eeaa
where $c_{n}^{\beta, p}$  is as in (\ref{conWishart}). It follows that
\beaa
P\Big(\frac{\lambda_{max}}{p} \leq x\Big) & \leq & c_{n}^{\beta, p}(p x)^{\beta n(n-1)/2}\cdot \Big(\int_0^{px}y^{\frac{\beta}{2}(p-n+1)-1}\cdot
e^{-y/2}\,dy\Big)^{n}\\
& \leq &
c_{n}^{\beta, p}(px)^{\beta n(n-1)/2}\cdot\Big(\frac{2}{(\beta -x)p-\beta n}(p x)^{\frac{\beta}{2}(p-n+1)} e^{-px/2}\Big)^n\\
& \leq & C^n\cdot c_{n}^{\beta, p}\cdot (p x)^{\beta np/2}\cdot p^{-n}\cdot e^{-npx/2}
\eeaa
as $n$ is sufficiently large, where $C$ is a constant not depending on $n$ and the second inequality follows from  (\ref{use}). Consequently,
\bea\lbl{american}
\limsup_{n\to\infty}\frac{1}{np}\log P\Big(\frac{\lambda_{max}}{p} \leq x\Big) \leq -\frac{x}{2}+\frac{\beta}{2}\log x +\limsup_{n\to\infty}\frac{1}{np}\log \Big(p^{\beta np/2}c_{n}^{\beta, p}\Big).
\eea
From (\ref{conWishart}) we have
\bea
& & \frac{1}{np}\log \Big(p^{\beta np/2}c_{n}^{\beta, p}\Big) \nonumber\\
%&= & \frac{\beta}{2}\log n + \frac{1}{np}\log c_{n}^{\beta, p}\\
& \leq & \frac{\beta}{2}\log p -\frac{\beta}{2}\log 2+ \frac{1}{np}\log \Big(\prod_{j=1}^n\Gamma(\frac{\beta}{2}(p-n+j))\Big)^{-1} +\frac{1}{np}\log \Big(\prod_{j=1}^n\Gamma(1+\frac{\beta}{2}j)\Big)^{-1} + o(1) \nonumber\\
& & \lbl{empty}
\eea
as $n\to\infty.$ It is known from the paragraph below (\ref{flawless}) that $\Gamma(1+x) \geq x^xe^{-Cx}\geq 1$ for $x\geq e^{C},$ where $C$ is an universal constant. Then,
\beaa
\log \prod_{j=1}^n\Gamma(1+\frac{\beta}{2}j) & \geq  &  \sum_{j=1}^n\big(\frac{\beta}{2}j\big) \log \big(\frac{\beta}{2}j\big) -C\sum_{j=1}^n\frac{\beta}{2}j + O(1)\\
 & = & \frac{\beta}{2} \sum_{j=2}^nj\log j+O(n^2)
\eeaa
as $n\to\infty.$ Notice $\sum_{j=2}^nj\log j\geq \int_1^n x\log x\,dx=(\frac{1}{2}x^2\log x-\frac{x^2}{4})|_{1}^n=\frac{1}{2}n^2\log n+o(n^2).$ Hence,
\beaa
\log \prod_{j=1}^n\Gamma(1+\frac{\beta}{2}j) \geq \frac{\beta}{4}n^2\log n+O(n^2)
\eeaa
as $n\to\infty.$ From this and (\ref{lemonade}) we get that  the sum of the third and fourth terms in (\ref{empty}) is bounded by
\beaa
& & \frac{1}{np}\Big(-\frac{\beta}{4}(\log p)n(2p-n+1)- \frac{\beta}{2}(\log \frac{\beta}{2}-1)np + \frac{\beta}{4}(\log \frac{\beta}{2})n^2+\frac{1}{2}n\log p \nonumber\\
& & \ \ \ \ \ \ \ \ \ \ \ \ \ \ \ \ \ \ \ \ \ \ \ \ \ \ \ \ \ \ \ \ \ \ \ \ \ \ \ \ \ \ \ \ \ \ \ \ \ \ \ \ \ \ \ \ \ \ \ \ \ \   -\frac{\beta}{4}n^2\log n\Big) +o(1) \\
&= & -\frac{\beta}{2}\log p - \Big(\frac{\beta}{4}\Big)\frac{n}{p}\log \frac{n}{p} -\frac{\beta}{2}(\log \frac{\beta}{2}-1) +o(1)\\
& = & -\frac{\beta}{2}\log p  -\frac{\beta}{2}(\log \frac{\beta}{2}-1) +o(1)
\eeaa
as $n\to\infty.$ Therefore,
\bea\lbl{soda}
\limsup_{n\to\infty}\frac{1}{np}\log \Big(p^{\beta np/2}c_{n}^{\beta, p}\Big) \leq -\frac{\beta}{2}\log \beta +\frac{\beta}{2}.
\eea
This joint with (\ref{american}) yields that
\beaa
\limsup_{n\to\infty}\frac{1}{np}\log P\Big(\frac{\lambda_{max}}{p} \leq x\Big) \leq -\Big(\frac{x-\beta}{2} -\frac{\beta}{2}\log \frac{x}{\beta}\Big) < 0
\eeaa
since $0<x<\beta.$ Therefore,
\bea\lbl{dove}
\limsup_{n\to\infty}\frac{1}{p}\log P\Big(\frac{\lambda_{max}}{p} \leq x\Big)=-\infty
\eea
for all $0<x<\beta.$ To prove (\ref{graffiti}), without loss of generality, we assume $F\subset [0, \infty)$ and $\beta\notin F.$ Since $F$ is closed, then either $F\subset [0, a]$, $F \subset [b, \infty)$ or $F\subset [0, a] \cup [b, \infty)$ for some constants $a\in F$ and $b\in F$ with $0<a<\beta< b.$ Thus (\ref{graffiti}) follows trivially from (\ref{compared1}) and (\ref{dove}).\\

\noindent{{\it The proof of (\ref{deer}).}} To prove (\ref{deer}), it is enough to show
\bea\lbl{award}
\liminf_{n\to\infty}\frac{1}{p}\log P\Big(\frac{\lambda_{max}}{p} \in G\Big) \geq - I(x)
\eea
for all $x\in G,$ where $G$ is an open subset of $\mathbb{R}.$ If $x<\beta,$ then (\ref{award}) holds automatically since $I(x)=\infty$. If $x=\beta,$ noticing $I(x)=0$ if and only if $x=\beta,$  we then know from (\ref{graffiti}) that $\frac{\lambda_{max}}{p}\to \beta$ in probability as $n\to\infty,$  thus $P(\frac{\lambda_{max}}{p} \in G)\to 1,$ hence (\ref{award}) is also true.

Now assume $x>\beta.$ Since $G$ is open, choose constants $r,a,b$ with $\beta < r <a<x<b$ and $(a, b)\subset G.$ Recall (\ref{stone}) and (\ref{holy}). Under the restriction that $pa<\lambda_1<pb$ and $\lambda_n<\cdots <\lambda_2 < pr,$ we know
\beaa
& & \prod_{i=2}^n|\lambda_1-\lambda_i|^{\beta} \geq \big(p(a-r)\big)^{\beta (n-1)}\ \ \mbox{and}\\
& & \int_{pa}^{pb}\lambda_1^{\frac{\beta}{2}(p-n+1)-1}e^{-\lambda_1/2}\,d\lambda_1 \geq p(b-a) (pa)^{\frac{\beta}{2}(p-n+1)-1}e^{-pb/2}.
\eeaa
Then, by the same argument as in (\ref{here}), we have
\bea
P\Big(\frac{\lambda_{max}}{p} \in G\Big) &\geq  & P\Big(a<\frac{\lambda_{1}}{p} <b,\,   \frac{\lambda_2}{p}   <r\Big)   \nonumber \\
% & = &  \int_{nx<\lambda_1,\,  \lambda_1>\cdots >\lambda_p>0}g(\lambda_1, \cdots, \lambda_p)\,d\lambda_1\cdots \,d\lambda_p \nonumber\\
& = & A_n \cdot \int_{pa}^{pb}B_n\,d\lambda_1 \cdot \int_{pr>\lambda_2>\cdots >\lambda_n>0}L_n(\lambda_2, \cdots, \lambda_n)\,d\lambda_2\cdots \,d\lambda_n\nonumber\\
& \geq &  A_n \cdot (p(a-r))^{\beta(n-1)}\cdot p(b-a)\cdot (pa)^{\frac{\beta}{2}(p-n+1)-1}e^{-pb/2}\cdot \tilde{P}\big(\frac{\lambda_2}{p}<r\big) \nonumber\\
& &  \lbl{pasture}
\eea
where $A_n$, $B_n$ and $L_n$ are defined in (\ref{stone}) and (\ref{holy}), and $\tilde{P}\big(\frac{\lambda_2}{p}<r\big)$ stands for the probability of $\{\frac{\lambda_2}{p}<r\}$ with the underlying probability distribution having density function $L_n(\lambda_2, \cdots, \lambda_n).$ Observe from (\ref{holy}) that the original beta-Laguerre ensemble with parameter $(n, p_n, \beta)$ becomes $L_n(\lambda_2, \cdots, \lambda_n)$ with parameter $(n-1, p_n-1, \beta).$ Since $\lim_{n\to\infty}p_n/n=\infty$, we know that $p_n':=p_{n+1}-1\to\infty$ and $p_n'/n\to\infty$ as $n\to\infty$, according to the arguments below (\ref{award}), we have $\lim_{n\to\infty}\tilde{P}\big(\frac{\lambda_2}{p}<r\big)=1$ since $r>\beta.$  Thus,
\bea\lbl{strange}
\frac{1}{p}\log P\Big(\frac{\lambda_{max}}{p} \in G\Big)\geq \frac{1}{p}\log A_n + \frac{\beta}{2}  \frac{n\log p}{p} +\frac{\beta}{2}\log p +\frac{\beta}{2}\log a - \frac{b}{2} + o(1)
\eea
as $n\to\infty.$ By (\ref{outrage}), the right hand side of the above is equal to
\beaa
-\frac{\beta}{2}(\log 2 -1 +\log \frac{\beta}{2})  + \frac{\beta}{2}\log a -\frac{1}{2}b + K_n
\eeaa
where
\beaa
K_n & = & \big(-\frac{\beta}{2}\frac{n\log n}{p} - \frac{\beta}{2}\log p\Big) + \Big(\frac{\beta}{2}\frac{n\log p}{p} + \frac{\beta}{2}\log p \Big)  +o(1)\\
&= & -\frac{\beta}{2}\frac{n}{p}\log \frac{n}{p} +o(1)\to 0
\eeaa
as $n\to\infty.$ Now taking $\liminf_{n\to\infty}$ for the both sides of the inequality in (\ref{strange}), and then letting $a \uparrow x$ and $b\downarrow x,$ we arrive at
\beaa
\liminf_{n\to\infty}\frac{1}{p}\log P\Big(\frac{\lambda_{max}}{p} \in G\Big) \geq -\frac{\beta}{2}(\log 2 -1 +\log \frac{\beta}{2})  + \frac{\beta}{2}\log x -\frac{1}{2}x=-I(x)
\eeaa
which gives (\ref{award}) for $x>\beta.$\ \ \ \ \ \ \ $\blacksquare$\\

\noindent\textbf{Proof of Theorem \ref{done1}}. It is easy to check that  $I(x)> 0$ for all $x\ne \beta$, $I(\beta)=0$, $\{I(x)\leq c\}$ is compact for any $c\geq 0,$ and $I(x)$ is strictly decreasing on $(0,\beta].$ Now, to prove the theorem, we need to show that
\bea
& & \limsup_{n\to\infty}\frac{1}{p}\log P\Big(\frac{\lambda_{min}}{p} \in F\Big) \leq -\inf_{x\in F} I(x)\lbl{graffiti1}\ \ \mbox{and}\\
& & \liminf_{n\to\infty}\frac{1}{p}\log P\Big(\frac{\lambda_{min}}{p} \in G\Big) \geq -\inf_{x\in G} I(x)\lbl{deer1}
\eea
for any closed set $F\subset \mathbb{R}$ and open set $G \subset \mathbb{R}.$

\noindent{{\it The proof of (\ref{graffiti1}).}} Obviously, the joint density function of the order statistics $\lambda_{min}=\lambda_{(n)}<\cdots <\lambda_{(1)}$ is $g(\lambda_1, \cdots, \lambda_n)=n!f_{n, \beta}(\lambda_{1}, \cdots, \lambda_n)$ for all $\lambda_1>\cdots >\lambda_n.$ Write
\bea\lbl{stone2}
& & g(\lambda_1, \cdots, \lambda_n)=\underbrace{\frac{n!c_{n}^{\beta, p}}{(n-1)!c_{n-1}^{\beta, p-1}}}_{A_n}\cdot \Big(\underbrace{\lambda_{n}^{\frac{\beta}{2}(p-n+1)-1}
e^{-\frac{1}{2}\lambda_{n}}\prod_{i=1}^{n-1}(\lambda_{i}-\lambda_n)^{\beta}}_{B_n}\Big)\cdot L_n(\lambda_1, \cdots, \lambda_{n-1})\nonumber\\
& &
\eea
where
\bea\lbl{holy2}
& & L_n(\lambda_1, \cdots, \lambda_{n-1})\nonumber\\
&= & (n-1)!c_{n-1}^{\beta, p-1}\prod_{1\leq i<j\leq
n-1}|\lambda_i-\lambda_j|^{\beta}\cdot\prod_{i=1}^{n-1}\lambda_i^{\frac{\beta}{2}((p-1)-(n-1)+1)-1}\cdot
e^{-\frac{1}{2}\sum_{i=1}^{n-1}\lambda_i}.
\eea
Observe that $\prod_{i=1}^{n-1}|\lambda_{n}-\lambda_i|^{\beta} \leq \lambda_1^{\beta (n-1)}$ for all $\lambda_1>\cdots >\lambda_n> 0.$ This gives
\beaa
B_n \leq (p M)^{\beta n}\cdot \lambda_n^{\frac{\beta}{2}(p-n+1)-1}e^{-\lambda_n/2}
\eeaa
provided $\lambda_1\leq pM \geq 1.$ By Theorem \ref{pack}, for any $M>\beta $,  we know that
\bea\lbl{royal}
P\Big(\frac{\lambda_{max}}{p} \geq M\Big) \leq e^{-p \gamma_M}
\eea
as $n$ is sufficiently large, where $\gamma_M=\frac{M-\beta}{4} -\frac{\beta}{4}\log \frac{M}{\beta}.$ Thus, for any $0<x<\beta,$ we have
\bea
  P\Big(\frac{\lambda_{min}}{p} \leq x\Big)& \leq & P\Big(\frac{\lambda_{min}}{p} \leq x, \frac{\lambda_{max}}{p}  \leq M) + e^{-p\gamma_M}\nonumber\\
& \leq  &2e^{-p\gamma_M} \vee  \int_{pM\geq   \lambda_1>\cdots >\lambda_n,\, \lambda_n\leq px}g(\lambda_1, \cdots, \lambda_n)\,d\lambda_1\cdots \,d\lambda_n \nonumber\\
& \leq & 2e^{-p\gamma_M} \vee (p M)^{\beta n}\cdot A_n \cdot \int_{0}^{p x}\lambda_n^{\frac{\beta}{2}(p-n+1)-1}e^{-\lambda_n/2}\,d\lambda_n \cdot\nonumber\\
 & & \ \ \ \ \ \ \ \ \ \ \ \ \ \ \ \ \ \ \ \  \ \ \ \ \ \ \ \ \ \ \ \  \int_{\lambda_1>\cdots >\lambda_{n-1}>0}L_n(\lambda_1, \cdots, \lambda_{n-1})\,d\lambda_1\cdots \,d\lambda_{n-1}\nonumber\\
& = & 2e^{-p\gamma_M} \vee \underbrace{(pM)^{\beta n}\cdot A_n \cdot \int_{0}^{p x}\lambda_n^{\frac{\beta}{2}(p-n+1)-1}e^{-\lambda_n/2}\,d\lambda_n}_{D_n}\lbl{way1}
\eea
since $L_n(\lambda_1, \cdots, \lambda_{n-1})$ is a probability density function.
By (\ref{use}), the last integral is bounded by $\frac{2}{(\beta-x)p-\beta n}(p x)^{\frac{\beta}{2}(p-n+1)} e^{-p x/2}$ as $n$ is sufficiently large. Thus, by (\ref{outrage}),
\beaa
\log D_n & \leq & \beta n\log (p M)-\frac{\beta}{2}p\big(\log 2 -1+\log \frac{\beta}{2}\big) -\frac{\beta}{2}n\log n -\frac{\beta}{2}p\log p \nonumber\\
& & \ \ \ \ \ \ \ \ \ \ \ +\log \frac{2}{ (\beta-x)p-\beta n}
+ \frac{\beta}{2}(p-n+1)(\log p + \log x)-\frac{p x}{2} + o(p)\nonumber\\
&= & \eta p -\frac{\beta}{2}\Big(\frac{n}{p}\log \frac{n}{p}\Big)p +o(p)
\eeaa
as $n\to\infty,$ where $\eta=-\frac{\beta}{2}\big(\log 2 -1+\log \frac{\beta}{2}\big) + \frac{\beta}{2}\log x-\frac{x}{2}.$  Since  $p\gg n$ it follows that
%Then, use the assumption {\red Start from here} $p=o(n/\log n)$ to have
\beaa
\limsup_{n\to\infty}\frac{1}{p}\log D_n & \leq & -\frac{\beta}{2}\big(\log 2 -1+\log \frac{\beta}{2}\big) + \frac{\beta}{2}\log x-\frac{x}{2}\\
& = & -\Big(\frac{x-\beta}{2} - \frac{\beta}{2}\log \frac{x}{\beta}\Big).
\eeaa
Consequently, by the notation $\gamma_M$ and (\ref{way1}) we get
%Note that $(p\log p)\vee(p\log n)=o(n)$ since $p=o(n/\log n).$ This joint with the above concludes that
\bea\lbl{compared2}
\limsup_{n\to\infty}\frac{1}{p}\log P\Big(\frac{\lambda_{min}}{p} \leq x\Big) \leq -\Big[\Big(\frac{x-\beta}{2}-\frac{\beta}{2}\log \frac{x}{\beta}\Big)\wedge  \Big(\frac{M-\beta}{4} -\frac{\beta}{4}\log \frac{M}{\beta}\Big)\Big]
\eea
for any $M>\beta$.  Send $M\to \infty$ to have
\bea\lbl{compared3}
\limsup_{n\to\infty}\frac{1}{p}\log P\Big(\frac{\lambda_{min}}{p} \leq x\Big) \leq - I(x)
%\Big(\frac{x-\beta}{2}-\frac{\beta}{2}\log \frac{x}{\beta}\Big)
\eea
for any $0<x<\beta.$

Now, if $\frac{1}{p}\lambda_{min}\geq x>\beta$ and $\frac{1}{p}\lambda_{max}\leq M$ for some $M>x$, by the definition of $f_{n,\beta}(\lambda)$ in (\ref{bWishart}), we see
\beaa
f_{n,\beta}(\lambda)\leq c_{n}^{\beta, p}(p M)^{\beta n(n-1)/2}\cdot \prod_{i=1}^n\lambda_i^{\frac{\beta}{2}(p-n+1)-1}\cdot
e^{-\frac{1}{2}\sum_{i=1}^n\lambda_i}
\eeaa
where $c_{n}^{\beta, p}$  is as in (\ref{conWishart}). It follows that
\beaa
P\Big(\frac{\lambda_{min}}{p} \geq x, \frac{\lambda_{max}}{p} \leq M\Big) & \leq & c_{n}^{\beta, p}(p M)^{\beta n(n-1)/2}\cdot \Big(\int_{px}^{\infty}y^{\frac{\beta}{2}(p-n+1)-1}\cdot
e^{-y/2}\,dy\Big)^{n}\\
& \leq &
c_{n}^{\beta, p}(pM)^{\beta n(n-1)/2}\cdot\Big((p x)^{\frac{\beta}{2}(p-n+1)} e^{-px/2}\Big)^n\\
& \leq & C^{n^2}\cdot c_{n}^{\beta, p}\cdot (p x)^{\beta np/2}\cdot e^{-npx/2}
\eeaa
as $n$ is sufficiently large, where $C=C(x, M)>0$ is a constant not depending on $n$, and the second inequality follows from  (\ref{seminar}) with $b-2\alpha-2=p(x-\beta)+\beta n -\beta -2\to\infty$. This implies
\bea\lbl{american1}
& & \limsup_{n\to\infty}\frac{1}{np}\log P\Big(\frac{\lambda_{min}}{p} \geq x, \frac{\lambda_{max}}{p} \leq M\Big) \nonumber\\
&\leq & -\frac{x}{2}+\frac{\beta}{2}\log x +\limsup_{n\to\infty}\frac{1}{np}\log \Big(p^{\beta np/2}c_{n}^{\beta, p}\Big).
\eea
This joint with (\ref{soda}) gives
\beaa
\limsup_{n\to\infty}\frac{1}{np}\log P\Big(\frac{\lambda_{min}}{p} \geq x, \frac{\lambda_{max}}{p} \leq M\Big) \leq -\Big(\frac{x-\beta}{2} -\frac{\beta}{2}\log \frac{x}{\beta}\Big) < 0
\eeaa
since $x>\beta.$ Therefore, by (\ref{royal}),
\beaa
& & \limsup_{n\to\infty}\frac{1}{p}\log P\Big(\frac{\lambda_{min}}{p} \geq x\Big)\nonumber\\
&\leq &  \limsup_{n\to\infty}\frac{1}{p}\Big[ \log P\Big(\frac{\lambda_{min}}{p} \geq x, \frac{\lambda_{max}}{p} \leq M\Big)\vee \log P\Big(\frac{\lambda_{max}}{p} \geq M\Big)\Big]\nonumber\\&=&-\infty\vee \Big(-\frac{M-\beta}{4} + \frac{\beta}{4}\log
 \frac{M}{\beta}\Big)
\eeaa
for all $M>x>\beta.$ By letting $M\to \infty$, we see that
\bea\lbl{dove1}
\limsup_{n\to\infty}\frac{1}{p}\log P\Big(\frac{\lambda_{min}}{p} \geq x \Big) =-\infty
\eea
for all $x>\beta.$ Since $P(\lambda_{min}>0)=1,$ by the same argument as in the paragraph below (\ref{dove}), we get (\ref{graffiti1}) from (\ref{compared3}) and (\ref{dove1}).

\noindent{{\it The proof of (\ref{deer1}).}} In order to prove (\ref{deer1}), we only need to show
\bea\lbl{award1}
\liminf_{n\to\infty}\frac{1}{p}\log P\Big(\frac{\lambda_{min}}{p} \in G\Big) \geq - I(x)
\eea
for all $x\in G \cap (0, \beta],$ where $G$ is an open subset of $\mathbb{R}.$ If $\beta\in G$, since $I(x)=0$ if and only if $x=\beta,$  we then know from (\ref{graffiti1}) that $\frac{\lambda_{min}}{p}\to \beta$ in probability as $n\to\infty,$  thus $P(\frac{\lambda_{min}}{p} \in G)\to 1,$ hence (\ref{award1}) is true for $x=\beta.$

Now it is enough to prove (\ref{award1}) for all $x\in G \cap (0, \beta) \ne \emptyset.$ For such $x,$ since $G$ is open, choose constants $a,b, r$ with $0< a<x<b<r<\beta$ and $(a, b)\subset G.$ Review (\ref{stone2}) and (\ref{holy2}). Under the restriction that $pa<\lambda_n<pb$ and $p r<\lambda_{n-1}<\cdots <\lambda_1 ,$ we have
\beaa
& & \prod_{i=1}^{n-1}|\lambda_i-\lambda_n|^{\beta} \geq \big(p(r-b)\big)^{\beta (n-1)}\ \ \mbox{and}\\
& & \int_{pa}^{pb}\lambda_n^{\frac{\beta}{2}(p-n+1)-1}e^{-\lambda_n/2}\,d\lambda_n \geq p(b-a) (pa)^{\frac{\beta}{2}(p-n+1)-1}e^{-pb/2}.
\eeaa
Then, by the same argument as in (\ref{pasture}), we have
\beaa
P\Big(\frac{\lambda_{min}}{p} \in G\Big) &\geq  & P\Big(a<\frac{\lambda_{n}}{p} <b,\,   \frac{\lambda_{n-1}}{p} >r\Big)    \nonumber\\
% & = &  \int_{nx<\lambda_1,\,  \lambda_1>\cdots >\lambda_p>0}g(\lambda_1, \cdots, \lambda_p)\,d\lambda_1\cdots \,d\lambda_p \nonumber\\
& = & A_n \cdot  \int_{ \lambda_1>\cdots >\lambda_{n-1}>pr}\big(\int_{pa}^{pb}B_n\,d\lambda_n\big) L_n(\lambda_1, \cdots, \lambda_{n-1})\,d\lambda_1\cdots \,d\lambda_{n-1}\nonumber\\
& \geq &  A_n \cdot (p(r-b))^{\beta(n-1)}\cdot p(b-a)\cdot (pa)^{\frac{\beta}{2}(p-n+1)-1}e^{-pb/2}\cdot \tilde{P}\big(\frac{\lambda_{n-1}}{p}>r\big)  \nonumber\\
& &
\eeaa
where $A_n$, $B_n$ and $L_n$ are defined in (\ref{stone2}) and (\ref{holy2}), and $\tilde{P}\big(\frac{\lambda_{n-1}}{p}>r\big)$ stands for the probability of $\{\frac{\lambda_{n-1}}{p}>r\}$ with the underlying probability distribution having density function $L_n(\lambda_1, \cdots, \lambda_{n-1}).$ Noticing $p_n':=p_{n+1}-1\to\infty$ and $p_n'/n\to\infty$ as $n\to\infty$ since $\lim_{n\to\infty}p_n/n=\infty.$ Then, by (\ref{graffiti1}) and the argument between (\ref{pasture}) and (\ref{strange}), we have  $\lim_{n\to\infty}\tilde{P}\big(\frac{\lambda_{n-1}}{p}>r\big)=1$ since $r<\beta.$  Thus,
\bea\lbl{strange1}
\frac{1}{p}\log P\Big(\frac{\lambda_{min}}{p} \in G\Big)\geq \frac{1}{p}\log A_n + \frac{\beta}{2}  \frac{n\log p}{p} +\frac{\beta}{2}\log p +\frac{\beta}{2}\log a - \frac{b}{2} + o(1)
\eea
as $n\to\infty.$  By (\ref{outrage}), the right hand side of the above is equal to
\beaa
-\frac{\beta}{2}(\log 2 -1 +\log \frac{\beta}{2})  + \frac{\beta}{2}\log a -\frac{1}{2}b + K_n
\eeaa
where
\beaa
K_n & = & \big(-\frac{\beta}{2}\frac{n\log n}{p} - \frac{\beta}{2}\log p\Big) + \Big(\frac{\beta}{2}\frac{n\log p}{p} + \frac{\beta}{2}\log p \Big)  +o(1)\\
&= & -\frac{\beta}{2}\frac{n}{p}\log \frac{n}{p} +o(1)\to 0
\eeaa
as $n\to\infty.$ Now taking $\liminf_{n\to\infty}$ for the both sides of (\ref{strange1}), and then letting $a \uparrow x$ and $b\downarrow x,$ we arrive at
\beaa
\liminf_{n\to\infty}\frac{1}{p}\log P\Big(\frac{\lambda_{min}}{p} \in G\Big) \geq -\frac{\beta}{2}(\log 2 -1 +\log \frac{\beta}{2})  + \frac{\beta}{2}\log x -\frac{1}{2}x=-I(x)
\eeaa
for all $x\in G \cap (0, \beta),$ which concludes (\ref{award1}).\ \ \ \ \ \ \ $\blacksquare$\\

\subsection{Proof of Theorem \ref{main1}}\lbl{ProofLDP}

\noindent To prove Theorem \ref{main1} next, we need to review some terminology. Let $\mathcal{M}(\mathbb{R})$ be the collection of the Borel probability measures defined on  $\mathbb{R}$ associated with the standard weak topology, that is, $\mu_n$ converges to $\mu$ weakly as $n\to\infty$ if and only if $\lim_{n\to\infty}\int f(x)\, \mu_n(dx) = \int f(x)\, \mu(dx)$  for every bounded and continuous function $f(x)$ defined on $\mathbb{R}$, where $\{\mu, \mu_n;\, n\geq 1\} \subset \mathcal{M}(\mathbb{R}).$ For further reference, see, e.g., chapter 11 from
Dudley (2002). When we mention open and closed sets in $\mathcal{M}(\mathbb{R})$ in the following, the corresponding topology is the weak topology.\\

\noindent\textbf{Proof of Theorem \ref{main1}}. By Theorem 1.3 from Ben Arous and Guionnet (1997), $I_{\beta}(\nu)$ is a good rate function, that is,   $I_{\beta}(\nu)\geq 0$ for all $\nu \in \mathcal{M}(\mathbb{R})$ and $\{\nu \in \mathcal{M}(\mathbb{R});\, I_{\beta}(\nu)\leq l\}$ is compact under the weak topology for any $l\geq 0.$ So we only need to show
\bea\lbl{disc}
\limsup_{n\to\infty}\frac{1}{n^2}\log P(\mu_n \in F) \leq -\inf_{\nu\in F} I_{\beta}(\nu)
\eea
for any closed set $F \subset \mathcal{M}(\mathbb{R})$ and
\bea\lbl{kid}
\liminf_{n\to\infty}\frac{1}{n^2}\log P(\mu_n \in G) \geq -\inf_{\nu\in G} I_{\beta}(\nu)
\eea
for any open set $G\subset \mathcal{M}(\mathbb{R}).$

\noindent{\it The proof of (\ref{disc})}. Define
\bea\lbl{shell}
E_n(\epsilon)=\Big\{\max_{1\leq i \leq n}|\frac{\lambda_i}{p}-\beta|< \sqrt{2\beta}\epsilon\Big\}
\eea
for $0<\epsilon <\frac{1}{4}(\sqrt{\beta} \wedge 1).$ By Theorems \ref{pack} and \ref{done1}, there exists a constant $\delta=\delta(\epsilon)>0$ such that
\bea\lbl{orange}
P(E_n(\epsilon)^c) \leq P\Big(\frac{\lambda_{max}}{p} \geq  \beta+\sqrt{2\beta}\epsilon\Big) + P\Big(\frac{\lambda_{min}}{p}\leq  \beta-\sqrt{2\beta}\epsilon\Big)  \leq  2e^{-p\delta}
\eea
as $n$ is sufficiently large. By (\ref{bWishart}),
\beaa
& & P(\mu_n \in F)\\
& \leq & 2e^{-p\delta} + P(\{\mu_n \in F\}\cap E_n(\epsilon))\\
& = & 2e^{-p\delta} + c_{n}^{\beta, p}\cdot \int_{\{\mu_n \in F\}\cap E_n(\epsilon)} \prod_{1\leq i<j\leq
n}|\lambda_i-\lambda_j|^{\beta}\cdot\prod_{i=1}^n\lambda_i^{\frac{\beta}{2}(p-n+1)-1}\cdot
e^{-\frac{1}{2}\sum_{i=1}^n\lambda_i}\,d\lambda_1\cdots d\lambda_n.
\eeaa
Since $x_i=\sqrt{\frac{p}{2\beta}}(\frac{\lambda_i}{p}-\beta).$ Then $\lambda_i=p(\beta +\sqrt{\frac{2\beta}{p}}x_i)$ with $x_i>-\sqrt{\beta p/2}$ for $1\leq i \leq n.$ It follows that
\bea
& & P(\mu_n \in F)\nonumber\\
& \leq & 2e^{-p\delta} + C_{n}^{\beta, p}\cdot \int_{\{\mu_n \in F\}\cap E_n(\epsilon)'} \prod_{1\leq i<j\leq
n}|x_i-x_j|^{\beta}\cdot\prod_{i=1}^n\Big(1+\sqrt{\frac{2}{p\beta}}x_i\Big)^{\frac{\beta}{2}(p-n+1)-1}\nonumber\\
& & \ \ \ \ \ \ \ \ \ \ \ \ \ \ \ \ \ \ \ \ \ \ \ \ \ \ \ \ \ \ \ \ \ \ \ \ \ \ \ \ \ \ \ \ \ \ \ \ \ \ \ \ \ \ \ \ \ \ \ \ \ \ \ \
\cdot e^{-\sqrt{\frac{\beta p}{2}}\sum_{i=1}^nx_i}\,dx_1\cdots dx_n\lbl{cotton}
\eea
where
\bea
& &
%\nu_n=\frac{1}{p}\sum_{i=1}^p\delta_{x_i/\sqrt{p}}\ \ \mbox{and}\ \
E_n(\epsilon)' =\Big\{\frac{\max_{1\leq i \leq n}|x_i|}{\sqrt{p}}< \epsilon\Big\}\ \ \mbox{and}\lbl{gram}\\
& & C_{n}^{\beta, p}=c_{n}^{\beta, p}\cdot (2\beta p)^{n(n-1)\beta/4}\cdot (p\beta)^{\frac{\beta}{2}n(p-n+1)-n}\cdot e^{-np\beta/2}\cdot (2\beta p)^{n/2}.\lbl{encounter}
\eea
By inequality $\log (1+x) \leq x -\frac{x^2}{2} +\frac{x^3}{3}$ for all $x>-1,$ we have
\beaa
& & \log \prod_{i=1}^n\Big(1+\sqrt{\frac{2}{p\beta}}x_i\Big)^{\frac{\beta}{2}(p-n+1)-1}\\
&\leq & \Big(\frac{\beta}{2}(p-n+1)-1\Big)\Big(\sqrt{\frac{2}{p\beta}}\sum_{i=1}^nx_i- \frac{1}{\beta p}\sum_{i=1}^nx_i^2 + \frac{2\sqrt{2}}{3\beta^{3/2}p^{3/2}}\sum_{i=1}^nx_i^3\Big)
\eeaa
for $x_i>-\sqrt{\beta p/2}$ with $i=1,\cdots, n.$ Now, on $E_n(\epsilon)'$ we have $\sum_{i=1}^n|x_i|\leq n\sqrt{p}\epsilon,$ it follows that

\bea\lbl{termU}
U_n:=\Big(\frac{\beta}{2}(p-n+1)-1\Big)\Big(\sqrt{\frac{2}{p\beta}}\sum_{i=1}^nx_i\Big)=\sqrt{\frac{\beta p}{2}}\sum_{i=1}^nx_i + \epsilon\cdot O(n^2)
\eea
as $n\to\infty.$ Similarly, on $E_n(\epsilon)'$ we have $\sum_{i=1}^n|x_i|^2\leq np\epsilon^2,$ which leads to
\bea\lbl{termV}
V_n:=-\Big(\frac{\beta}{2}(p-n+1)-1\Big)\cdot  \frac{1}{\beta p}\sum_{i=1}^nx_i^2=-\frac{1}{2}\sum_{i=1}^nx_i^2 +\epsilon\cdot O(n^2)
\eea
as $n\to\infty$ since $0<\epsilon<1.$ By the same argument, on $E_n(\epsilon)'$ we have $\sum_{i=1}^n|x_i|^3\leq \sqrt{p}\epsilon\sum_{i=1}^n|x_i|^2,$ hence,
\bea\lbl{termW}
W_n:=\Big(\frac{\beta}{2}(p-n+1)-1\Big)\Big( \frac{2\sqrt{2}}{3\beta^{3/2}p^{3/2}}\sum_{i=1}^n|x_i|^3\Big)\leq \frac{\epsilon}{\sqrt{\beta}}\sum_{i=1}^n x_i^2.
\eea
Combining all the above we get
\beaa
\log \prod_{i=1}^n\Big(1+\frac{x_i}{\beta\sqrt{p}}\Big)^{\frac{\beta}{2}(p-n+1)-1}  \leq  U_n+V_n +W_n \leq \sqrt{\frac{\beta p}{2}}\sum_{i=1}^n x_i-\frac{\alpha_{\epsilon}}{2}\sum_{i=1}^nx_i^2 + \epsilon\cdot O(n^2)
\eeaa
on $E_n(\epsilon)'$ as $n\to\infty,$ where
\bea\lbl{writes}
\alpha_{\epsilon}:=1-\frac{2\epsilon}{\sqrt{\beta}}>0
\eea
for $0<\epsilon <\frac{1}{4}(\sqrt{\beta} \wedge 1).$ From (\ref{cotton}) we see that
\bea
& & P(\mu_n \in F)\nonumber\\
&\leq & 2e^{-p\delta} + C_{n}^{\beta, p}\cdot\exp\{\epsilon\cdot O(n^2)\}\cdot\int_{\mu_n \in F}\prod_{1\leq i<j\leq
n}|x_i-x_j|^{\beta}\cdot e^{-\frac{\alpha_{\epsilon}}{2}\sum_{i=1}^nx_i^2}\,dx_1\cdots dx_n \ \ \ \ \ \ \ \ \lbl{flips}
%\nonumber\\
%& = & 2e^{-n\delta} + C_{n}^{\beta, p}\cdot\exp\{\epsilon\cdot O(p^2)\}\cdot \alpha_{\epsilon}^{-\frac{\beta}{2}p(p-1)-p}\nonumber\\
%& & \ \ \ \ \ \ \ \ \ \ \ \ \ \ \ \ \ \  \cdot \int_{\mu_{n,1} \in  F}\prod_{1\leq i<j\leq
%p}|y_i-y_j|^{\beta}\cdot e^{-\frac{1}{2}\sum_{i=1}^py_i^2}\,dy_1\cdots dy_p\lbl{trade}
\eea
%by setting $y_i=\alpha_{\epsilon}x_i$ for $1\leq i \leq p,$ where
%\bea\lbl{immediate}
%\mu_{n,1}=\frac{1}{p}\sum_{i=1}^p\delta_{y_i/(\alpha_{\epsilon}\sqrt{p})}.
%\eea
as $n$ is sufficiently large. Let $\lambda_1, \cdots, \lambda_n$ have the density function $f_{\beta}(\lambda_1, \cdots, \lambda_n)$ as in (\ref{betaHermite}). Set  $y_i=\lambda_i/\sqrt{\alpha_{\epsilon}}$ for $i=1,\cdots, n.$ We know that $(y_1, \cdots, y_n)$ has density
\bea\lbl{gesture}
h_{\beta}(y_1, \cdots, y_n):=\alpha_{\epsilon}^{\frac{\beta}{4}n(n-1) + \frac{n}{2}}K_n^{\beta}\prod_{1\leq i <  j \leq n}|y_i-y_j|^{\beta}\cdot e^{-\frac{\alpha_{\epsilon}}{2}\sum_{i=1}^ny_i^2}
\eea
for $(y_1,\cdots, y_n)\in \mathbb{R}^n$. By Corollary 5.1 from  Ben Arous and Guionnet (1997) (taking $a=\alpha_{\epsilon}$, $V(x)=\alpha_{\epsilon}x^2$ and $f(x)\equiv 0$) that $\nu_n=\frac{1}{n}\sum_{i=1}^n\delta_{y_i/\sqrt{n}}$
satisfies the LDP with speed $\{n^2;\, n\geq 1\}$ and rate function $\alpha_{\epsilon}\cdot I_{\beta/\alpha_{\epsilon}}(\nu)$
where
%$I_{b}(\nu)$ is defined in (\ref{LDPrate}).
\bea\lbl{LDPrate}
I_{b}(\nu)=\frac{1}{2}\int_{\mathbb{R}^2}g_{b}(x,y)\,\nu(dx)\,\nu(dy) +\frac{b}{4}\log \frac{b}{2}-\frac{3}{8}b
\eea
for any $b>0$ and
\bea\lbl{elm}
g_{b}(x,y)=
\begin{cases}
\frac{1}{2}(x^2+y^2) - b \log |x-y|, & \text{if $x\ne y$;}\\
+\infty, & \text{otherwise.}
\end{cases}
\eea
We see from (\ref{flips}) that
\bea
& &  P(\mu_n \in F)\nonumber\\
&\leq & 2e^{-p\delta} + \frac{C_{n}^{\beta, p}}{K_n^{\beta}}\cdot\exp\{\epsilon\cdot O(n^2)\}\cdot \alpha_{\epsilon}^{-\frac{\beta}{4}n(n-1)-\frac{n}{2}}\cdot P(\nu_{n} \in  F) \lbl{wave}
\eea
as $n\to\infty$. It follows that
\bea\lbl{laugh}
& & \limsup_{n\to\infty}\frac{1}{n^2}\log P(\mu_n \in F)\nonumber\\
&\leq & \limsup_{n\to\infty}\frac{1}{n^2}\log \Big(\frac{C_{n}^{\beta, p}}{K_n^{\beta}}\Big) -\inf_{\nu\in F}\Big\{\alpha_{\epsilon}\cdot I_{\beta/\alpha_{\epsilon}}(\nu)\Big\} -\frac{\beta}{4}\log \alpha_{\epsilon} + \epsilon\cdot O(1)
\eea
where the condition $p \gg n^2$ is used in the  inequality. For the constant and rate function above, we have the following facts.
\begin{lemma}\lbl{ask} If  $p=p_n \gg n$, then $\log \frac{C_n^{\beta,p}}{K_n^{\beta}} = O(n+n^3p^{-1})$ as $n\to\infty.$ In particular, it is of order $O(n)$ if $p\gg n^2$.
%If  $p=p_n \gg n^2$, then $\log \frac{C_n^{\beta,p}}{K_n^{\beta}} = O(n)$ as $n\to\infty.$
\end{lemma}

\begin{lemma}\lbl{almost} Let $I_b(\nu)$ be defined as in (\ref{LDPrate}). Let $A$ be a set of Borel probability measures on $\mathbb{R}.$ Set $J_{s}(A)=\inf_{\nu\in A} \{s\cdot I_{b/s}(\nu)\}$ for all $s>0$. If $\{s_n> 0;\, n\geq 1\}$ is a sequence with  $\lim_{n\to\infty}s_n=s\in (0, \infty),$ then $\lim_{n\to\infty}J_{s_n}(A)=J_{s}(A).$
\end{lemma}
The proofs of the two lemmas will be given at the end of this section. Let's continue now. From (\ref{laugh}) and Lemma \ref{ask},
\bea
& & \limsup_{n\to\infty}\frac{1}{n^2}\log P(\mu_n \in F)\nonumber\\
&\leq &  -\inf_{\nu\in F}\Big\{\alpha_{\epsilon}\cdot I_{\beta/\alpha_{\epsilon}}(\nu)\Big\} -\frac{\beta}{4}\log \alpha_{\epsilon} + \epsilon\cdot O(1)  \lbl{brown}
\eea
for any $\epsilon \in (0,\frac{1}{4}(\sqrt{\beta} \wedge 1 )).$ Now passing $\epsilon \downarrow 0,$ we know from (\ref{writes}) that $\alpha_{\epsilon} \to 1$, hence it follows from Lemma \ref{almost} that
\bea
 \limsup_{n\to\infty}\frac{1}{n^2}\log P(\mu_n \in F) \leq -\inf_{\nu \in F}I_{\beta}(\nu).\lbl{white}
%\nonumber\\
%&\leq & -\frac{\beta\log (2\beta)}{4} - \frac{1}{2\beta}\cdot\inf_{\nu\in F}I_{2\beta^2}(\nu) +\frac{\beta}{4}\log (2\beta)\nonumber\\
%& = & -\frac{1}{2\beta}\cdot\inf_{\nu\in F}I_{2\beta^2}(\nu)=-\inf_{\nu\in F}J(\nu)
%\lbl{white}
\eea
%with $J(\nu) = \frac{1}{2\beta}I_{2\beta^2}(\nu).$

\noindent{\it The proof of (\ref{kid})}. First, by the Taylor expansion, there exists $\epsilon_0 \in (0,1)$ such that
\bea\lbl{trivial}
\log (1+x)\geq x-\frac{1}{2}x^2-|x|^3
\eea
for all $|x|<\epsilon_0.$
%Define
%\beaa
%H_n(\epsilon)=\Big\{\max_{1\leq i \leq p}|\frac{\lambda_i}{n}-\beta|< \epsilon\beta\Big\}
%\eeaa
%for any $0<\epsilon<\epsilon_0.$
Without loss of generality, assume
\bea\lbl{ponding}
0< \epsilon_0 < \frac{1}{4} (\sqrt{\beta} \wedge 1).
\eea
Review (\ref{shell}) and (\ref{cotton}) to have
\bea
& & P(\mu_n \in G)\\
& \geq &  P(\{\mu_n \in G\}\cap E_n(\epsilon)) \nonumber\\
& = & c_{n}^{\beta, p}\cdot \int_{\{\mu_n \in G\}\cap E_n(\epsilon)} \prod_{1\leq i<j\leq
n}|\lambda_i-\lambda_j|^{\beta}\cdot\prod_{i=1}^n\lambda_i^{\frac{\beta}{2}(p-n+1)-1}\cdot
e^{-\frac{1}{2}\sum_{i=1}^n\lambda_i}\,d\lambda_1\cdots d\lambda_n  \nonumber\\
& = & C_{n}^{\beta, p}\cdot \int_{\{\mu_n \in G\}\cap E_n(\epsilon)'} \prod_{1\leq i<j\leq
n}|x_i-x_j|^{\beta}\cdot\prod_{i=1}^n\Big(1+\sqrt{\frac{2}{\beta p}}x_i\Big)^{\frac{\beta}{2}(p-n+1)-1}\nonumber\\
& & \ \ \ \ \ \ \ \ \ \ \ \ \ \ \ \ \ \ \ \ \ \ \ \ \ \ \ \ \ \ \ \ \ \ \ \ \ \ \ \ \ \ \ \ \ \ \ \ \ \ \ \ \ \ \ \ \ \ \ \ \ \ \ \
\cdot e^{-\sqrt{\frac{\beta p}{2}}\sum_{i=1}^nx_i}\,dx_1\cdots dx_n\lbl{cloth}
\eea
where $E_n(\epsilon)'$ is defined in (\ref{gram}) and $C_{n}^{\beta, p}$ in (\ref{encounter}).
%,  and
%\beaa
%H_n(\epsilon)' =\Big\{\frac{\max_{1\leq i \leq p}|x_i|}{\sqrt{n}}< \epsilon\beta\Big\}.
%\eeaa
Now, by the inequality in (\ref{trivial}), on $E_n(\epsilon)'$,  we have
\beaa
& & \log \prod_{i=1}^n\Big(1+\frac{x_i}{\beta\sqrt{p}}\Big)^{\frac{\beta}{2}(p-n+1)-1}\\
&\geq & \Big(\frac{\beta}{2}(p-n+1)-1\Big)\Big(\sqrt{\frac{2}{\beta p}}\sum_{i=1}^nx_i- \frac{1}{\beta p}\sum_{i=1}^nx_i^2 - \frac{2\sqrt{2}}{\beta^3p^{3/2}}\sum_{i=1}^n|x_i|^3\Big)\\
& = & U_n+V_n -3W_n
\eeaa
where $U_n, V_n$ and $W_n$ are defined in (\ref{termU}), (\ref{termV}) and (\ref{termW}), respectively. Thus,
%Comparing $H_n(\epsilon)'$ to  $E_n(\epsilon)'$ as in (\ref{gram}), replacing $\epsilon$ by $\beta \epsilon$ in (\ref{termU}), (\ref{termV}) and %(\ref{termW}), we obtain from (\ref{ponding}) that
\beaa
\log \prod_{i=1}^n\Big(1+\frac{x_i}{\beta\sqrt{p}}\Big)^{\frac{\beta}{2}(p-n+1)-1} \geq \sqrt{\frac{\beta p}{2}}\sum_{i=1}^nx_i-\frac{\gamma_{\epsilon}}{2}\sum_{i=1}^nx_i^2 + \epsilon\cdot O(n^2)
\eeaa
on $E_n(\epsilon)'$ as $n\to\infty,$ where
\bea\lbl{speaks}
\gamma_{\epsilon}:=1+\frac{6\epsilon}{\sqrt{\beta}}>0
\eea
for any $0<\epsilon<\epsilon_0.$ This joint with (\ref{cloth}) yields that
\beaa
& & P(\mu_n \in G)\nonumber\\
& \geq & C_{n}^{\beta, p}\cdot e^{\epsilon\cdot O(n^2)}\cdot  \int_{\{\mu_n \in G\}\cap E_n(\epsilon)'} \prod_{1\leq i<j\leq
n}|x_i-x_j|^{\beta}\cdot e^{-\frac{\gamma_{\epsilon}}{2}\sum_{i=1}^nx_i^2}\,dx_1\cdots dx_n\\
& = & \frac{C_{n}^{\beta, p}}{K_n^{\beta}}\cdot e^{\epsilon\cdot O(n^2)}\cdot \gamma_{\epsilon}^{-\frac{\beta}{4}n(n-1)-\frac{n}{2}}\cdot P_1(\{\mu_{n} \in  G\} \cap E_n(\epsilon)')
\eeaa
by the same arguments as  those in (\ref{gesture}) and (\ref{wave}), where $P_1$ stands for the probability such that $x=(x_1, \cdots, x_n)$ appearing in the definitions of $\mu_n$ and $E_n(\epsilon)'$ has the probability density function
\beaa
\tilde{h}_{\beta}(x_1, \cdots, x_n):=\gamma_{\epsilon}^{\frac{\beta}{4}n(n-1) + \frac{n}{2}}K_n^{\beta}\prod_{1\leq i <  j \leq n}|x_i-x_j|^{\beta}\cdot e^{-\frac{\gamma_{\epsilon}}{2}\sum_{i=1}^nx_i^2}
\eeaa
for $(x_1, \cdots, x_n)\in \mathbb{R}^n.$ From Lemma \ref{ask}, we get
\bea
& & \liminf_{n\to\infty}\frac{1}{n^2}\log P(\mu_n \in G) \nonumber\\
& \geq & \liminf_{n\to\infty}\frac{1}{n^2}\log P_1(\{\mu_{n} \in  G\} \cap E_n(\epsilon)') + \epsilon \cdot  O(1)-\frac{\beta}{4}\log \gamma_{\epsilon} \lbl{blue}
\eea
for any $0<\epsilon<\epsilon_0.$ Denote $J_{\epsilon}(G):=\inf_{\nu\in G} \{\gamma_{\epsilon}\cdot I_{\beta/\gamma_{\epsilon}}(\nu)\}.$ We claim that $(\ref{blue})$ implies
\bea\lbl{quiet}
\liminf_{n\to\infty}\frac{1}{n^2}\log P(\mu_n \in G)\geq -J_{\epsilon}(G)-\epsilon  + \epsilon \cdot  O(1)-\frac{\beta}{4}\log \gamma_{\epsilon}
\eea
for all $0<\epsilon<\epsilon_0.$ In fact, the above  is trivially true if $J_{\epsilon}(G)=\infty.$ Assume now $J_{\epsilon}(G)<\infty.$ Then, by the LDP discussion between (\ref{gesture}) and (\ref{wave}), we know $P_1(\mu_{n} \in  G)\geq \exp\{-n^2(J_{\epsilon}(G)+\epsilon)\}$ as $n$ is sufficiently large. From Lemma \ref{tennis},
\beaa
P_1(E_n(\epsilon)'^c) =P_1\Big(\max_{1\leq i \leq n}|x_i|\geq \sqrt{p}\epsilon\Big)=P\Big(\max_{1\leq i \leq n}|\lambda_i|\geq \sqrt{n}(\sqrt{p/n}\,\gamma_{\epsilon}^{1/2}\epsilon)\Big) \leq e^{-p\gamma_{\epsilon}\epsilon^2 /3}
\eeaa
as $n$ is sufficiently large since $p \gg n$, where $(\lambda_1, \cdots, \lambda_n):=\gamma_{\epsilon}^{1/2}\cdot (x_1, \cdots, x_n)$ has the joint probability density function $f_{\beta}(\lambda_1, \cdots, \lambda_n)$ as in (\ref{betaHermite}). Hence,
\beaa
P_1(\{\mu_{n} \in  G\} \cap E_n(\epsilon)') & \geq & P_1(\mu_{n} \in  G) - P_1(E_n(\epsilon)'^c)\geq e^{-n^2(J_{\epsilon}(G)+\epsilon)} - e^{-p\gamma_{\epsilon}\epsilon^2 /3}\\
& = &  e^{-n^2(J_{\epsilon}(G)+\epsilon)}(1+o(1))
\eeaa
as $n\to\infty$ by the condition $p \gg n^2.$ This and (\ref{blue}) lead to (\ref{quiet}). Finally, letting $\epsilon \downarrow 0$ in (\ref{quiet}), we have   (\ref{kid}) from Lemma \ref{almost}.\ \ \ \ \ \ \ \ $\blacksquare$\\

\noindent\textbf{Proof of Proposition \ref{LSD}.}  Let $\nu_{\beta}$ be the measure with density $g_{\beta}(x)=(\beta \pi)^{-1}\sqrt{2\beta -x^2}$ for any $|x|\leq \sqrt{2\beta}$,  and define $F_{\beta, \epsilon}= \{\mu\in \mathbb{R}; \, \rho(\mu, \nu_{\beta})\geq \epsilon\}$, where $\rho(\cdot, \cdot)$ is the Prohorov distance, see, e.g., chapter 11 from
Dudley (2002). Then, $F_{\beta, \epsilon}$ is a closed set under the weak topology. Recalling the definition of $\alpha_{\epsilon}$ in (\ref{writes}) and Lemma \ref{ask}, we have from (\ref{wave}) that, for $\epsilon$ small enough,
\beaa
P(\mu_n\in F_{\beta, \epsilon})& \leq & 2e^{-p\delta}+ Ce^{\epsilon O(n^2)}\cdot P(\nu_n\in F) \nonumber \\
& \leq &  2e^{-p\delta}+ Ce^{\epsilon O(n^2)}\cdot \exp\Big\{-(n^2/2)\cdot \inf_{\nu\in F}\{\alpha_{\epsilon}\cdot I_{\beta/\alpha_{\epsilon}}(\nu)\}\Big\} \nonumber
%\inf_{\nu \in F_{\beta,\epsilon}}I_{\beta}(\nu)}.
\eeaa
by the large deviation principle mentioned in (\ref{LDPrate}), where $\nu_n$ is defined between (\ref{gesture}) and (\ref{LDPrate}). By Lemma \ref{almost}, the infimum goes to a positive constant since $\nu_{\beta}\notin F$. So the desired result follows from the Borel-Cantelli lemma.\ \ \ \ \ $\blacksquare$\\

\noindent\textbf{Proof of Lemma \ref{ask}}. Review (\ref{encounter}).  Notice
\bea
\frac{C_n^{\beta,p}}{K_n^{\beta}}
&= &\frac{c_n^{\beta,p}}{K_n^{\beta}}\cdot (2\beta p)^{n(n-1)\beta/4}\cdot (p\beta)^{\frac{\beta}{2}n(p-n+1)-n}\cdot e^{-np\beta/2}\cdot (2\beta  p)^{n/2} \nonumber\\
& = & \underbrace{2^{-\frac{\beta}{2}np}\cdot(2\pi)^{n/2}\cdot\beta^{\frac{\beta}{2}n(p-n+1)-n}\cdot p^{\beta n(2p-n+1)/4-(n/2)}e^{-\beta np/2}\cdot (2\beta)^{\beta n(n-1)/4 + (n/2)}}_{D_n} \nonumber\\
& & \ \ \ \ \ \ \ \ \ \ \ \ \ \ \ \ \ \ \ \ \ \ \ \ \ \ \ \ \ \ \ \ \ \ \ \ \ \ \ \ \ \ \  \cdot \Big(\prod_{j=1}^n\Gamma(\frac{\beta}{2}(p-n+j))\Big)^{-1}.\lbl{sandwitch}
\eea
Observe
\beaa
\log D_n & = & -\frac{\beta}{2}(\log 2)np + \frac{n}{2}\log (2\pi) + \frac{\beta}{2}(\log \beta)n(p-n+1) -n\log \beta \\
& & \ \ \ \ \ \ \ \ \ \ \ \ \ \ \ \ \    + \frac{\beta}{4}n(2p-n+1)\log p -\frac{n}{2}\log p-\frac{\beta}{2} np + \big(\frac{\beta}{4}n(n-1)+\frac{n}{2}\big)\log (2\beta)\\
% & = & \frac{\beta}{4}(\log n)p(2n-p+1)  + \frac{\beta}{2}(\log \beta)p(n-p+1)  -\frac{\beta}{2} np + O(p)\\
& = & \frac{\beta}{4}(\log p)n(2p-n+1)  + \frac{\beta}{2}(\log \beta)n(p-n) -\frac{\beta}{2}(1+\log 2) np -\frac{n}{2}\log p\\
& &   \ \ \ \ \ \ \ \ \ \ \ \ \ \ \ \ \  \ \ \ \  \ \ \ \ \ \ \ \ \ \ \ \ \ \ \ \ \  \ \ \ \ \ \ \ \ \ \ \ \ \ \ \ \ \ +\frac{\beta\log (2\beta)}{4}n^2+ O(n)\\
& = & \frac{\beta}{4}(\log p)n(2p-n+1) -\frac{\beta}{4}(\log \frac{\beta}{2})n^2 + \frac{\beta}{2}(\log \beta-1-\log 2)np -\frac{n}{2}\log p+ O(n)
\eeaa
as $n\to\infty,$ where $O(n)$ is a function of $n$ which does not depend on $p$.  Joint this and (\ref{lemonade})  with  (\ref{sandwitch}) we conclude
\beaa
\log \frac{C_n^{\beta,p}}{K_n^{\beta}}
=   O\big(n+\frac{n^3}{p}\big)
\eeaa
as $p\gg n\to\infty$. In particular, it is of order $O(n)$ if $p\gg n^2.$\ \ \ \ \ \ \ \ $\blacksquare$\\

\noindent\textbf{Proof of Lemma \ref{almost}}. Recalling the expression of $I_b(\nu)$ in (\ref{LDPrate}) and (\ref{elm}), to prove the lemma, it is enough to show that
\bea\lbl{scan}
& & H_n:=\inf_{\nu\in A}\int_{\mathbb{R}^2}\big(a_n(x^2+y^2)- \log |x-y|\big)\, v(dx) v(dy)\nonumber\\
& \to & H:=\inf_{\nu\in A}\int_{\mathbb{R}^2}\big(a(x^2+y^2)- \log |x-y|\big)\, v(dx) v(dy)
\eea
as $n\to\infty,$ where $\{a_n>0;\, n\geq 1\}$ is a sequence with $\lim_{n\to\infty}a_n=a\in (0, \infty).$

We first claim that
%there exists a finite constant $C=C_t$ such that
\bea\lbl{document}
t(x^2+y^2)- \log |x-y| + C \geq \frac{t}{2}(x^2+y^2)
\eea
for  any $t>0$, $x>0$ and $y>0,$ where $C=C_t=-\frac{1}{2}\cdot \inf_{u>0} \{tu-\log u\}$ is finite. In fact,
\beaa
\frac{t}{2}(x^2+ y^2)-\log |x-y| \geq \frac{t}{2}(x-y)^2 - \frac{1}{2}\log (x-y)^2 \geq -C,
\eeaa
which yields (\ref{document}).
%where $\psi(u)=tu-\log u$ for $u>0$ and hence $C$ is a finite number.
The inequality in (\ref{document}) implies that
\bea\lbl{wallet}
& & \int_{\mathbb{R}^2}\big(t(x^2+y^2)- \log |x-y|\big)\, v(dx) v(dy)\ \ \mbox{is finite if and only if}\nonumber\\
& &  \ \int_{\mathbb{R}}x^2v(dx)<\infty\ \ \mbox{and}\ \ \int_{\mathbb{R}^2}\log |x-y|\, v(dx) v(dy)\  \ \mbox{is finite}.\ \ \ \ \ \ \
\eea
The inequality in (\ref{document}) also says that the sequence $\{H, H_n;\, n\geq 1\}$ is bounded below.

Now, if $H<\infty,$ take any $\nu\in A$ such that $\int_{\mathbb{R}^2}\big(a(x^2+y^2)- \log |x-y|\big)\, v(dx) v(dy)<\infty.$ By (\ref{wallet}),
\beaa
H_n\leq a_n\int_{\mathbb{R}^2}(x^2+y^2)v(dx) v(dy)- \int_{\mathbb{R}^2}\log |x-y|\, v(dx) v(dy)<\infty
\eeaa
for all $n\geq 1.$ Passing $n\to\infty$ and taking the infimum over all $\nu \in A$ for both sides, we get
\bea\lbl{hurry}
\limsup_{n\to\infty}H_n\leq H.
\eea
Obviously, the above is also true if $H=\infty.$

On the other hand, if $H<\infty,$ by (\ref{wallet}), $H_n<\infty$ for all $n\geq 1.$ For any $n\geq 1,$ take $\nu_n\in A$ with
\bea
H_n+\frac{1}{n} & \geq & \int_{\mathbb{R}^2}\big(a_n(x^2+y^2)- \log |x-y|\big) v_n(dx) v_n(dy)  \lbl{windy}\\
& \geq & \int_{\mathbb{R}^2}\big(2\delta(x^2+y^2)- \log |x-y|\big) v_n(dx) v_n(dy)\nonumber\\
& \geq & \delta\int_{\mathbb{R}^2}(x^2+y^2) v_n(dx) v_n(dy) -C  \nonumber
\eea
where $2\delta:=\inf\{a_n;\,n\geq 1\}$ and the last inequality follows from (\ref{document}) with $C=C_{\delta}$. This and (\ref{hurry}) imply that
\beaa
M:=\sup_{n\geq 1}\int_{\mathbb{R}^2}(x^2+y^2) v_n(dx) v_n(dy) < \infty.
\eeaa
Therefore, from (\ref{windy}),
\beaa
H_n + \frac{1}{n} & \geq & \int_{\mathbb{R}^2}\big(a(x^2+y^2)- \log |x-y|\big) v_n(dx) v_n(dy) -M |a_n-a|\\
& \geq & H - M |a_n-a|.
\eeaa
Letting $n\to\infty,$ we have $\liminf_{n\to\infty}H_n\geq H.$ Further, (\ref{wallet}) says that $H=\infty$ if and only if $H_n=\infty$ for all $n\geq 1.$ Thus, $\liminf_{n\to\infty}H_n\geq H$ also holds if $H=\infty.$ These and (\ref{hurry}) prove (\ref{scan}).\ \ \ \ \ \  $\blacksquare$

\section{Appendix}\lbl{check}
In this part, we verify the validity of (\ref{oversize}).

\noindent\textbf{Proof of (\ref{oversize})}. Set
\beaa
& & A=-\frac{\beta np}{2}\log 2 +\big(\frac{1}{2}\log \frac{\beta}{2}\big)n +\frac{1}{2}n\log p - \frac{1}{2}(\frac{\beta}{2}\log \frac{\beta}{2}-\frac{\beta}{2})n(2p-n+1)   \\
& & \ \ \ \ \ \ \ \ \ \ \ \ \ \ \ \ \ \ \ \ \ \ \ \ \ \ \ \ \ \ \ \ \ \ \ \ \ \ \ \ -\frac{\beta}{4}(\log p)n(2p-n+1)+\frac{\beta}{4}n(n-1);\\
& & B=  \big(\frac{1}{4}n(n-1)\beta +\frac{1}{2}n\big)\log (2\beta p)-\frac{1}{2}np\beta + \big(\frac{1}{2}n(p-n+1)\beta - n\big)\log (p\beta ).
\eeaa
To show (\ref{oversize}), by (\ref{gull}) and (\ref{sound}), we need to check that $A+B=0.$ First, separate the terms with $\log p$ from others in the expression of $B$ and then sort out for $\log \beta$ to have
\beaa
B &= & \big(\frac{1}{4}n(n-1)\beta +\frac{1}{2}n\big)\log p +\big(\frac{1}{4}n(n-1)\beta +\frac{1}{2}n\big)\log (2\beta) + \\
& & \big(\frac{1}{2}n(p-n+1)\beta - n\big)\log p +\big(\frac{1}{2}n(p-n+1)\beta - n\big)\log \beta -\frac{1}{2}np\beta\\
& = & \Big(\frac{\beta}{2}pn\log p-\big(\frac{1}{4}n(n-1)\beta +\frac{1}{2}n\big)\log p\Big) +B'
\eeaa
where
\beaa
B'=\frac{ (\log \beta) n}{4}\big((2p-n+1)\beta-2\big)   + \big(\frac{1}{4}n(n-1)\beta +\frac{1}{2}n\big)\log 2-\frac{1}{2}np\beta.
\eeaa
Comparing the two coefficients of $\log p$ in $A$ and $B$, we find that their sum is identical to $0$. Thus, by setting
\beaa
A'=-\frac{\beta np}{2}\log 2 +\big(\frac{1}{2}\log \frac{\beta}{2}\big)n  - \frac{1}{2}(\frac{\beta}{2}\log \frac{\beta}{2}-\frac{\beta}{2})n(2p-n+1)  +\frac{\beta}{4}n(n-1),
\eeaa
we only need to check $A'+B'=0.$ In fact, write
\bea
A' &= & -\frac{\beta np}{2}\log 2 +\big(\frac{1}{2}\log \frac{\beta}{2}\big)n  - \frac{(\log \beta) n}{4}(2p-n+1)\beta  + \frac{\beta(\log 2) n}{4}(2p-n+1) \nonumber\\
& & \ \ \ \ \ \ \ \ \ \ \ \ \ \ \ \ \ \ \ \ \ \ \ \ \ \ \ \ \ \ \ \ \ \   +\frac{\beta}{4}n(2p-n+1) +\frac{\beta}{4}n(n-1) \lbl{line}\\
& = &  - \frac{(\log \beta) n}{4}(2p-n+1)\beta + \frac{1}{2}(\log \beta)n + \Big(-\frac{\beta n(n-1)}{4}-\frac{n}{2}\Big)\log 2+\frac{1}{2}np\beta \nonumber
\eea
where in the second equality we first merge all terms with $\log 2$ together by  writing $\big(\frac{1}{2}\log \frac{\beta}{2}\big)n=\frac{1}{2}(\log \beta)n-\frac{1}{2}(\log 2)n,$ and then sum the last two terms in (\ref{line}) to obtain  $\frac{1}{2}np\beta.$ Now it is evident that $A'=-B'.$\ \ \ \ \ \ \ \ \ $\blacksquare$\\\\

\noindent {\bf Acknowledgments.} Dr. Alice Guionnet told the open problem to the first author on the large deviations for the eigenvalues of the Wishart matrices when the sample size is much smaller than the dimension of data at the ``{\it France-China Summer School and Conference: Random Matrix Theory and High-dimensional
Statistics, Changchun, China, 2011}". The problem is studied in Theorem \ref{main1}. We thank her very much. We also thank the associate editor and referees for their precious comments.

\baselineskip 10pt
\def\ref{\par\noindent\hangindent 30pt}

\end{document}